\documentclass[10pt]{amsart}
\usepackage{amssymb}
\usepackage{amscd}
\usepackage{comment}
\theoremstyle{plain}
\newtheorem{theorem}{Theorem}
\theoremstyle{definition}
\newtheorem{question}{Question}
\newtheorem*{notation}{Notation}
\newtheorem{proposition}{Proposition}[section]

\newtheorem{corollary}[proposition]{Corollary}
\newtheorem{lemma}[proposition]{Lemma}
\newtheorem{definition}[proposition]{Definition}
\newtheorem*{mainthm}{Main Theorem}
\theoremstyle{remark}
\newtheorem{claim}{Claim}

\DeclareMathOperator{\drop}{Drop}
\DeclareMathOperator{\wsat}{Wsat}
\DeclareMathOperator{\indec}{Indec}
\DeclareMathOperator{\Gap}{Gap}
\DeclareMathOperator{\refl}{Refl}
\DeclareMathOperator{\fil}{Fill}

\DeclareMathOperator{\nacc}{nacc}
\DeclareMathOperator{\Sk}{Sk}

\DeclareMathOperator{\otp}{otp}
\DeclareMathOperator{\Ch}{Ch}

\DeclareMathOperator{\comp}{Comp}

\DeclareMathOperator{\pr}{Pr}
\DeclareMathOperator{\Fill}{Fill}

\DeclareMathOperator{\ran}{ran}

\DeclareMathOperator{\cf}{cf}

\DeclareMathOperator{\ord}{Ord}
\newcommand{\sk}{\vskip.05in}
\DeclareMathOperator{\id}{id}

\newcommand{\restr}{\upharpoonright}

\newcommand{\uf}{\mathcal{U}}

\newcommand{\subs}{\subseteq}

\DeclareMathOperator{\acc}{acc}
\numberwithin{equation}{section}
\begin{document}
\title{Club-guessing, stationary reflection, and coloring theorems}
\author{Todd Eisworth}
\dedicatory{For Isabelle}
\address{Department of Mathematics\\
         Ohio University\\
         Athens, OH 45701}
\email{eisworth@math.ohiou.edu}
 \keywords{square-brackets partition relations, minimal walks, successor of singular cardinal, club-guessing, stationary reflection}
 \subjclass{03E02}
\thanks{The author acknowledges support from NSF grant DMS 0506063.}
\date{\today}
\begin{abstract}
We obtain strong coloring theorems at successors of singular cardinals from failures of certain instances
of simultaneous reflection of stationary sets. Along the way, we establish new results in club-guessing and in the general theory of ideals.
\end{abstract}

\maketitle

\section{Introduction}

Our main results establish that certain failures of simultaneous reflection for stationary subsets of the successor of a singular cardinal can be used to generate strong versions of negative square-brackets partition relations. This represents a substantial improvement over previous results --- we get stronger conclusions from much weaker hypotheses. We are able to obtain this by way of progress in two areas. First, our primary advance allows us to obtain a version of the main result of~\cite{nsbpr} which holds even for successors of singular cardinals of countable cofinality. We accomplish this by changing the sort of club-guessing sequence used in those papers, and by refining our arguments so that they work in this new context. On another front, we extend results of Shelah connecting the structure of ideals with simultaneous reflection of stationary sets, and then combine this with the other advance to deduce the main theorem.

In this introductory section, our goal is to provide enough background so that we are both able to state our main result precisely and situate it in context.  We assume the reader has an acquaintance with basic set-theoretic notation; any of the standard references (say \cite{jechbook},\cite{kunenbook}, or \cite{akibook}) are more than sufficient.  Some of the objects of importance to us have been denoted in a variety of ways throughout the literature, so we take a moment to set our conventions.

\begin{notation}
Let $\lambda$ and $\kappa$ be cardinals.
\begin{enumerate}
\item $[\lambda]^{\kappa} :=\{A\subs\lambda: |A|=\kappa\}$.
\sk
\item $S^\lambda_\kappa:=\{\alpha<\lambda:\cf(\alpha)=\kappa\}$.
\end{enumerate}
We also utilize minor variants of the above notation, but these should all be self-explanatory.
\end{notation}

One more convention is quite important for us: by {\em ``$I$ is an ideal on $\kappa$''}, we mean  {\bf ``$I$ is a proper ideal on~$\kappa$ containing all bounded subsets of~$\kappa$''}. We have chosen this notation for convenience, since it eliminates many trivialities.

 We now move on to the background material; we begin with the {\em square-brackets partition relation} $\kappa\rightarrow [\lambda]^\mu_\theta$ of Erd\H{o}s, Hajnal, and Rado~\cite{EHR}, which means that for any function $F:[\kappa]^{\mu}\rightarrow\theta$,
(which we refer to as a {\em coloring}) we can find a set $H\subs\kappa$ of cardinality $\lambda$ for which
\begin{equation*}
\ran(F\restr [H]^\mu)\subsetneqq\theta,
\end{equation*}
that is, when we restrict the function $F$ to $[H]^\mu$, at least one color is omitted.  These square-brackets partition relations arise naturally when one investigates the extent to which Ramsey's Theorem generalizes to uncountable sets, as the statement ``$\kappa\rightarrow[\lambda]^\mu_\theta$'' asserts that a very weak form of Ramsey's Theorem holds at the cardinal~$\kappa$.

The negations of these partition relations are quite strong combinatorial hypotheses in their own right, and it is for this reason that they have been studied extensively by set theorists:
 \begin{quote}
 \small{
 The moral is that if one knows an ordinary negative partition relation then it is often worth asking whether a stronger assertion, a negative square bracket relation, is also true.} \cite{ehmr}
 \end{quote}
Consider for example, the statement $\kappa\nrightarrow[\kappa]^2_\kappa$.  This
asserts the existence of a function $F:[\kappa]^2\rightarrow \kappa$ (a {\em coloring} of the pairs from $\kappa$ using $\kappa$ colors) with the property that $F$ assumes every possible value on (the pairs from) every unbounded subset of $\kappa$. This says that Ramsey's Theorem fails in a very spectacular way at~$\kappa$: it is possible to color  the pairs of ordinals less than $\kappa$ utilizing $\kappa$ colors such that any subset of $\kappa$ of size $\kappa$ is {\em completely inhomogeneous} with respect to the coloring.  Theorems that assert the existence of ``complicated'' colorings of pairs (or other finite sets) of ordinals (where the exact meaning of ``complicated'' depends on context) can be conveniently grouped under the sobriquet ``coloring theorems'', and our main theorem is one such.

In Appendix~1 of~\cite{cardarith}, Shelah systematically studies several combinatorial principles stronger than the negated square-brackets relations discussed above, and our results are most naturally stated using his notation.

\begin{definition}
\label{prdef}
If $\lambda$ is an infinite cardinal, and $\kappa+\theta\leq\mu\leq\lambda$, then $\pr_1(\lambda,\mu,\kappa,\theta)$ asserts the existence of function $c:[\lambda]^2\rightarrow\kappa$ such that whenever we are given a collection $\langle t_\alpha:\alpha<\mu\rangle$ of pairwise disjoint elements of $[\lambda]^{<\theta}$ as well as an ordinal $\varsigma<\kappa$,  then there are $\alpha<\beta$ for which $c\restr t_\alpha\times t_\beta$ is constant with value~$\varsigma$.
\end{definition}

Notice that for $\theta>2$, the relation $\pr_1(\lambda,\mu,\kappa,\theta)$ implies $\lambda\nrightarrow [\mu]^2_\kappa$, as one may take the sets $t_\alpha$ to be singletons; this justifies our referring to theorems that establish instances of $\pr_1$ as coloring theorems.

We will break off our discussion of these relations for a moment to pick up another thread that is important for our work, namely, reflection of stationary sets.  Jech's~\cite{jechhandbook} gives a nice introduction to this topic, and the paper~\cite{cfm} of Cummings, Foreman, and Magidor is also an excellent resource.

\begin{definition}
Let $S$ be a stationary subset of an uncountable regular cardinal~$\kappa$.
\begin{enumerate}
\item We say {\em $S$ reflects at $\alpha$} if $\alpha<\kappa$ has uncountable cofinality and $S\cap\alpha$ is stationary in~$\alpha$.
    \sk
\item $\refl(S)$ holds if every stationary subset of $S$ reflects at some~$\alpha$.
\sk
\item $S$ is non-reflecting if $S$ does not reflect at any~$\alpha$.
\end{enumerate}
\end{definition}

We observe here that any cardinal of the form $\kappa^+$ for regular $\kappa$ has a non-reflecting stationary subset: consider the set $S^{\kappa^+}_\kappa$ of all ordinals below $\kappa^+$ of cofinality~$\kappa$. This set does not reflect, as any $\alpha<\kappa$ of uncountable cofinality contains a closed unbounded subset of ordinals each of cofinality less than $\kappa$.  The situation at successors of singular cardinals is much more delicate, and we will have much to say about this later.

It has been known for a long time that there are connections between stationary reflection and coloring theorems.  For example, Tryba~\cite{tryba} (and independently Hugh Woodin) established that $\kappa\nrightarrow [\kappa]^{<\omega}_\kappa$ whenever $\kappa$ has a non-reflecting stationary subset. Stevo Todor{\v{c}}evi{\'c}~\cite{acta}  was able to use his technique of {\em minimal walks} to improve this to colorings of pairs:  $\kappa\nrightarrow [\kappa]^2_\kappa$ whenever $\kappa$ has a non-reflecting stationary set, and hence, in particular, $\kappa^+\nrightarrow [\kappa^+]^2_{\kappa^+}$ whenever $\kappa$ is a regular cardinal. Shelah~\cite{186} simplified Todor{\v{c}}evi{\'c}'s argument a bit and obtained some other generalizations.

Research on coloring theorems for successors of regular cardinals has continued. For example, Shelah~\cite{572} has established that $\pr_1(\kappa^{+2},\kappa^{+2},\kappa^{+2},\kappa)$ holds for every regular cardinal~$\kappa$, while Justin Moore~\cite{lspace}  obtained a significant strengthening of Todor{\v{c}}evi{\'c}'s result for $\aleph_1$ and solved a long-standing open problem in general topology.

Given our understanding of successors of regular cardinals, the following question is natural:

\begin{question}
To what extent are analogous results true for successors of singular cardinals?
\end{question}

This question is still very mysterious, in large part due to the fact that the combinatorics of successors of singular
cardinals is very sensitive to underlying assumptions in set theory.  For example, if $\mu$ is singular and $\mu^+$ has a non-reflecting stationary subset, then the work of Todor{\v{c}}evi{\'c} can be brought to bear and we can conclude $\mu^+\nrightarrow[\mu^+]^2_{\mu^+}$. However,  assuming the existence of large cardinals, Magidor~\cite{magidor} established the consistency of $\refl(\mu^+)$ for $\mu$, and so the situation at successors of singular cardinals differs greatly from that at successors of regular cardinals.

Back in 1978, Shelah~\cite{shelahjonsson} was able to establish $\aleph_{\omega+1}\nrightarrow[\aleph_{\omega+1}]^2_{\aleph_{\omega+1}}$ from the assumption that $2^{\aleph_0}<\aleph_\omega$.  His later development of pcf theory let him eliminate the additional assumption (see Chapter~II of~\cite{cardarith}), while also extending the class of cardinals for which such results hold.  Morever,  in {\sf ZFC}, Shelah (Conclusion~4.1 in Chapter~II of~\cite{cardarith}) was able to prove
\begin{equation}
\pr_1(\mu^+,\mu^+,\cf(\mu),\cf(\mu))\text{ holds for any singular cardinal }\mu.
\end{equation}
Notice that this is a strong negative partition relation using $\cf(\mu)$ colors; the question of whether
 or not this can be improved to a coloring using~$\mu^+$ colors (that is, whether $\pr_1(\mu^+,\mu^+,\mu^+,\cf(\mu))$ holds for~$\mu$ singular) is still very much open, as are the related questions concerning $\mu^+\nrightarrow[\mu^+]^2_{\mu^+}$ and $\mu^+\nrightarrow[\mu^+]^{<\omega}_{\mu^+}$ . The main result of this paper obtains  strong colorings under very weak assumptions, but to describe the precise situation, we need some notation from~\cite{cfm}.

 \begin{definition}
 Let $\kappa$ be a regular uncountable cardinal, and let $S$ be a stationary subset of $\kappa$.
 Then $\refl(\theta, S)$ holds if for every sequence $\langle T_i:i<\theta\rangle$ of stationary
 subsets of $S$, there exists an $\alpha<\kappa$ such that each $T_i$ reflects at~$\alpha$.
Similarly, we say $\refl(<\theta, S)$ holds if for every $\sigma<\theta$ and every sequence $\langle T_i:i<\sigma\rangle$
 of stationary subsets of~$S$, there exists an $\alpha<\kappa$ such that each~$T_i$ reflects at~$\alpha$.
 \end{definition}

 We can now state our main theorem, a theorem illustrating the connection between failures of simultaneous reflection of stationary sets and the existence of complicated colorings.

 \begin{mainthm}
{\em Assume $\mu$ is a singular cardinal.
\begin{enumerate}
\sk
\item If $\refl(<\cf(\mu),S^{\mu^+}_{\geq\theta})$ fails for some $\theta<\mu$, then $\pr_1(\mu^+,\mu^+,\theta,\cf(\mu))$ holds.

    \sk

\item If $\refl(<\!\cf(\mu), S^{\mu^+}_{\geq\theta})$ fails for arbitrarily large $\theta<\mu$, then we obtain both $\pr_1(\mu^+,\mu^+,\mu, \cf(\mu))$  and $\mu^+\nrightarrow[\mu^+]^2_{\mu^+}$.

\sk

\item If $\cf(\mu)>\aleph_0$ and $\refl(<\cf(\mu), S^{\mu^+}_{\geq\theta})$ fails for arbitrarily large $\theta<\mu$, then~(2) can be improved to $\pr_1(\mu^+,\mu^+,\mu^+,\cf(\mu))$.
\sk

\item If $\cf(\mu)=\aleph_0$ and $\refl(<\cf(\mu), S^{\mu^+}_{\geq\theta})$ fails for arbitrarily large $\theta<\mu$, then there is a function $d:[\mu^+]^3\rightarrow\mu^+$ such that whenever $\langle t_\alpha:\alpha<\mu^+\rangle$ is a pairwise disjoint family of finite subsets of $\mu^+$ and $\varsigma<\mu^+$, there are $\alpha<\beta<\gamma$ such that
 \begin{equation*}
(\forall\epsilon\in t_\alpha)(\forall \zeta\in t_\beta)(\forall \xi\in t_\gamma)[ d(\epsilon,\zeta,\xi)=\varsigma].
\end{equation*}
\end{enumerate}
}
\end{mainthm}

 The proof builds on work of the author and Shelah~\cite{535, 819, nsbpr}.  In particular, the main theorem of~\cite{nsbpr}, when taken together with results in sections 2 and 3 of the current paper, is strong enough to imply
 our main theorem in the case that~$\mu$ has uncountable cofinality.  The paper~\cite{819} partially extended
 the results of~\cite{nsbpr} to the case where $\mu$ has countable cofinality (the combination of the main theorem in~\cite{819} with the results of sections 2 and 3 in the current paper results in something weaker than $\pr_1(\mu^+,\mu^+,\mu^+,\cf(\mu))$); the current paper gets around the obstacle in~\cite{819} caused by the countable cofinality of~$\mu$.

The proof of most parts of the main theorem will proceed according to the following sketch:

\medskip
\noindent{\sf Sketch of argument}

\medskip
\noindent Let $\lambda=\mu^+$ for $\mu$ singular.
\begin{enumerate}
\item In {\sf ZFC} we prove that there are a function $c:[\lambda]^2\rightarrow \lambda$ and an ideal $I$ (related to club-guessing) such that
whenever $\langle t_\alpha:\alpha<\lambda\rangle$ is a family of pairwise disjoint elements of $[\lambda]^{<\cf(\mu)}$, for
almost every (in the sense of $I$) $\beta^*<\lambda$ there are $\alpha<\beta$ such that $c\restr t_\alpha\times t_\beta$
is constant with value $\beta^*$.
\sk
\item If the desired colorings fail to exist, then we can conclude ideal $I$ from~(1) possesses some strong combinatorial properties.
\sk
\item These properties are strong enough to imply the needed instances of simultaneous reflection of stationary subsets of~$\lambda$.
\sk
\end{enumerate}

\medskip

The main theorem of~\cite{nsbpr} establishes~(1) for the case where the cofinality of~$\mu$ is uncountable, while
the main theorem of~\cite{819} gives a conclusion weaker than~(1). Statement~(2) will be obtained from~(1) by modifying
one of the arguments from~\cite{nsbpr}. Statement~(3) relies on some extending some results from section~3 of Chapter~IV of Shelah's~\cite{cardarith}; the results seem to be new, though they have a ``folklore-ish'' flavor.

Moving on, we now give an account of the organization of this paper.
In sections 2 and 3, we prove a series of results culminating in a proof of~(3) above. Section~4 proves a club-guessing
result needed to define the ideal~$I$. Section~5 pins down the assumptions we need for our main theorems, and lays
some groundwork for later arguments as well. Section~6 introduces much of the minimal walks machinery we need, as well as providing
a proof of a crucial preliminary result.  Section~7 furnishes some background material concerning scales and elementary
submodels, and Section~8 provides a proof of~(1), as well as a few other odds and ends we use in the proof of our main theorem.  Section~9 closes the paper by finishing the proof of our primary result.

Finally, some historical remarks are in order. Clearly, the research in this paper rests on that of Saharon Shelah,
 and although the main theorem is an advance in our knowledge, the proof is obtained via a synthesis of techniques drawn
 from several places in his vast body of work. The line of investigation which resulted in this paper originated in our work with
 Shelah in~\cite{535}, where I noticed that the argument in Section~4 of Chapter~III in~\cite{cardarith} which purported to extend the main coloring theorem
 in that section to cardinals of the form $\mu^+$ for $\aleph_0=\cf(\mu)<\mu$ did not work.  A close study of the problem led
 to the isolation of the coloring theorems in~\cite{nsbpr}, and in the joint paper~\cite{819} we exploited a combinatorial
 trick to provide a partial rescue of the result from~\cite{cardarith}.  Some years earlier, Shelah had suggested the
 idea of ``off-center'' club-guessing as a possible way of repairing the error in~\cite{cardarith}, but nothing came of the idea
 at the time because the combinatorics did not work as we hoped.  While writing up~\cite{819}, the author realized that
 the ``off-center'' approach might be viable when combined with the combinatorics of that paper. The work in sections~3 through~9
 shows that this was indeed the case --- in Theorem~\ref{mainfilterthm} we get a version of the main result from~\cite{nsbpr} which holds for successors of singular cardinals of countable cofinality.  The jump from this to getting results on simultaneous reflection of stationary sets is based on other work of Shelah --- Section 2 and 3 of this paper can be viewed as teasing out additional consequences from arguments appearing in Section~3 of Chapter~III from~\cite{cardarith}.

\section{Weak saturation and indecomposability}

Our goal in this section is a modest one --- we will take two fairly well-known properties of ideals and show, using elementary arguments, that their conjunction is equivalent to several useful properties. The results appearing in this section have a ``folklore flavor'' and may be known; Shelah maps out some of the implications appearing here in Chapter~III of his~\cite{cardarith}, but he stops short of proving that many of the properties he was considering are in fact equivalent.  Theorem~\ref{thm1} at the end of this section gathers the results in a single place.

Before proving anything, we take a moment to discuss some notational conventions and expressions that will appear over and over in our results. This terminology is fairly standard, so readers familiar with such things can just skip ahead.

\begin{definition}
Let $I$ be an ideal on $\kappa$.
\begin{enumerate}
\item $I^*$ denotes the filter dual to $I$.
\sk
\item Expressions of the form ``$\varphi$ holds for $I$-almost all $\alpha<\kappa$'' mean that the set of $\alpha<\kappa$ for which $\varphi$ holds is in $I^*$.  If $I$ is clear from context, we may omit explicit reference to it and say only ``$\varphi$ holds for almost all $\alpha<\kappa$''.
\sk
\item If $A$ and $B$ are subsets of $\kappa$, then $A\subs_I B$ ({\em $A$ is a subset of $B$ modulo $I$}) means that $A\setminus B\in I$, and we say $A=_I B$ ({\em $A$ and $B$ are equal modulo~$I$}) if both $A\subs_I B$ and $B\subs_I A$.
\sk
\item Similarly, if $f$ and $g$ are ordinal-valued functions with domain~$\kappa$, we say that {\em $f\leq_I g$}
if $f(\alpha)\leq g(\alpha)$ for almost all $\alpha<\kappa$. The expression  $f=_I g$ is also given the obvious meaning.
\end{enumerate}
\sk
Variants of the above notation should also be interpreted in the canonical fashion.
\end{definition}

We look now at the first of two properties which form our main interest in this section.

\begin{definition}
Let $I$ be an ideal on the cardinal $\kappa$, and let $\theta$ be a cardinal. The ideal $I$ is {\em weakly $\theta$-saturated} if there is no partition of $\kappa$ into $\theta$ disjoint $I$-positive sets.
\end{definition}

Weak saturation provides a measure of how closely the filter $I^*$ dual to $I$ comes to being an ultrafilter. For example, $I^*$ is an ultrafilter if and only if $I$ is weakly $2$-saturated, while
the ideal of bounded subsets of $\kappa$ fails to be weakly $\kappa$-saturated. This concept has not been studied as systematically
as its better-known relative ``saturation'', but we shall see that it is quite important in its own right.
The following  observation starts us on our way.

\begin{proposition}
\label{ob1}
Let $I$ be an ideal on the cardinal $\kappa$. The following statements are equivalent for a cardinal $\theta$:
\begin{enumerate}
\sk
\item $I$ is weakly $\theta$-saturated.
\sk
\item Any $\subs$-increasing $\theta$-sequence $\langle A_i:i<\theta\rangle$ of subsets of $\kappa$ is eventually constant modulo $I$.
\sk
\item If $\theta\leq\cf(\tau)$, then any $\subs$-increasing $\tau$-sequence $\langle A_i:i<\tau\rangle$ of
 subsets of $\kappa$ is eventually constant modulo $I$.
\end{enumerate}
\end{proposition}
\begin{proof}
The proof is trivial. For example, one proves that (1) implies (2) via contradiction: If (2) fails, then we
can refine $\langle A_i:i<\theta\rangle$ to a subsequence $\langle B_i:i<\theta\rangle$ with the property that
$B_i\subs_I B_{i+1}$ and $B_{i+1}\setminus B_i\notin I$ for each $i<\theta$.  The collection $\langle B_{i+1}\setminus B_i:i<\theta\rangle$ quickly leads to a contradiction of~(1).
\end{proof}

The next concept we need for our discussion is that of {\em indecomposability}, defined as follows:

\begin{definition}
Let $I$ be an ideal on $\kappa$, and let $\theta$ be a cardinal. The ideal $I$ is said to be {\em $\theta$-indecomposable} if whenever
$\langle A_i:i<\theta\rangle$ is a $\theta$-sequence of subsets of $\kappa$ with $\bigcup_{i<\theta}A_i\notin I$, there is a set $w\subs\theta$
of cardinality less than $\theta$ with $\bigcup_{i\in w}A_i\notin I$. In the case where $\theta$ is a regular cardinal, this is equivalent to the statement
that the ideal $I$ is closed under {\em increasing} union of length $\theta$.
\end{definition}

Indecomposability has been considered many times in the literature (see \cite{indec1} and \cite{indec2} for example), although one usually finds the definition phrased in terms of filters rather than ideals, and very often the authors restrict themselves to considering ultrafilters rather than the more general case.

\begin{proposition}
\label{obs2}
Let $I$ be an ideal on the cardinal $\kappa$.  The following two statements are equivalent for a regular cardinal $\theta$:
\begin{enumerate}
\sk
\item $I$ is $\theta$-decomposable

\sk

\item There is a $\theta$-sequence $\langle A_i: i<\theta\rangle$ of $I$-equivalent $I$-positive sets with
\begin{equation}
\label{eqn1}
\bigcap_{i<\theta}\bigcup_{i\leq j<\theta}A_i  = \emptyset.
\end{equation}
\sk
\end{enumerate}
\end{proposition}

\begin{proof}
Suppose $I$ is $\theta$-decomposable. Since $\theta$ is regular, this means there is an {\em increasing} sequence $\langle B_i:i<\theta\rangle$
of elements of $I$ whose union $B$ is $I$-positive.  If we define
\begin{equation*}
A_i:= B\setminus B_i,
\end{equation*}
then $\langle A_i:i<\theta\rangle$ has the required properties.

For the other direction, suppose we are given $\langle A_i:i<\theta\rangle$ as in (2).  Let $B=A_0$, and define
\begin{equation*}
B_i = B\setminus \bigcup_{i\leq j<\theta}A_j.
\end{equation*}
Note that $B_i\in I$ as $B$ and $A_i$ are equivalent modulo $I$.  The sequence $\langle B_i:i<\theta\rangle$ is also increasing,
and by~(\ref{eqn1}) we have
\begin{equation*}
B=\bigcup_{i<\theta}B_i.
\end{equation*}
Since $B\notin I$, we conclude that $I$ is $\theta$-decomposable.
\end{proof}

Our main concern in this section is a consideration of the conjunction of ``$\theta$-indecomposable'' and ``weakly $\theta$-saturated'' for a regular cardinal~$\theta$. The next proposition shows us that this combination
is has some strength --- it is equivalent to improvements of both Proposition~\ref{ob1} and Proposition~\ref{obs2}.

\begin{proposition}
\label{prop1}
 Let $I$ be an ideal on the cardinal $\kappa$. The
following statements are equivalent for a regular cardinal~$\theta$.
\begin{enumerate}
\sk
\item $I$ is weakly $\theta$-saturated and $\theta$-indecomposable.

\sk

\item Whenever $\langle B_i:i<\theta\rangle$ is an increasing $\theta$-sequence of subsets of $\kappa$, there is an $i^*$ such that
\begin{equation}
i^*\leq i<\theta\Longrightarrow B_i =_I \bigcup_{j<\theta}B_j.
\end{equation}

\sk

\item Whenever $\langle A_i:i<\theta\rangle$ is a $\theta$-sequence of $I$-positive subsets of~$\kappa$, we have
\begin{equation}
\label{eqn2}
\bigcap_{i<\theta}\bigcup_{i\leq j<\theta} A_j \neq \emptyset.
\end{equation}
\end{enumerate}
\end{proposition}
\begin{proof}
Assume (1) holds, and let $\langle B_i:i<\theta\rangle$ be an increasing $\theta$-sequence of subsets of~$\kappa$.  Since $I$
is weakly $\theta$-saturated, we know that this sequence is eventually constant modulo $I$, so fix $i^*<\theta$ with the property that
\begin{equation*}
i^*\leq i<\theta\Longrightarrow B_{i^*} = B_i \mod I.
\end{equation*}
Let $B=\bigcup_{i<\theta}B_i$, and note that (2) is established if we can prove that $B\setminus B_{i^*}$ is in $I$. This is done
through $\theta$-indecomposability, for $B\setminus B_{i^*}$ can be expressed as the union of an increasing $\theta$-sequence
of elements of $I$:
\begin{equation*}
B\setminus B_{i^*} = (\bigcup_{i^*\leq i<\theta}B_i)\setminus B_{i^*} =\bigcup_{i^*\leq i<\theta}(B_i\setminus B_{i^*}).
\end{equation*}

Next, assume that condition (2) holds for our ideal, and assume by way of contradiction that
 $\langle A_i:i<\theta\rangle$ is a $\theta$-sequence of $I$-positive subsets of~$\kappa$ for which~(\ref{eqn2}) fails.
For each $i<\theta$, define
\begin{equation*}
B_i:=\kappa\setminus \bigcup_{i\leq j<\theta}A_j.
\end{equation*}
The sequence $\langle B_i:i<\theta\rangle$ is increasing, and furthermore $\bigcup_{i<\theta}B_i = \kappa$ because~(\ref{eqn2}) fails.
By (2), there is an $i^*$ such that $B_{i^*}=_I \kappa$.
This is a contradiction, as $A_{i^*}$ is $I$-positive and disjoint to $B_{i^*}$.

To finish the proof, we show that the failure of (1) implies the failure of (3).  This is easily done --- if $I$ is not weakly $\theta$-saturated
than any partition of $\kappa$ into $\theta$ disjoint $I$-positive sets will contradict (3), and if $I$ is $\theta$-decomposable
then we contradict~(3) by way of Proposition~\ref{obs2}.
\end{proof}

Notice that condition~(2) says something stronger than the conclusion of Proposition~\ref{ob1} ---  the
sequence of sets is not only eventually constant modulo~$I$, it is the case that eventually the individual sets $B_i$
 are equal to the union of the entire sequence modulo~$I$.

Condition~(3) says that any ``point--$<\theta$'' collection of $I$-positive sets has size less than $\theta$, for
given a collection of $\theta$ (or more) $I$-positive sets, there is a subcollection of size $\theta$ with non-empty intersection.

Since prime ideals are weakly $2$-saturated, we obtain the following (known) characterization
 of $\theta$-indecomposability in the context of ultrafilters.

\begin{corollary}
Let $\uf$ be an ultrafilter on some cardinal $\kappa$, and let $\theta$ be a regular cardinal.
Then $\uf$ is $\theta$-decomposable if and only if it is $(\theta,\theta)$-regular.
\end{corollary}
\begin{proof}
By definition $\uf$ is $(\theta,\theta)$-regular if and only if there is a family $\{A_i:i<\theta\}$ of elements of $\uf$ with
the property that the intersection of any subfamily of size $\theta$ is empty.  This is precisely the condition~(\ref{eqn2}), and
so the result follows immediately as the ideal dual to $\uf$ is trivially weakly $\theta$-saturated.
\end{proof}

Our next move takes us to the realm of functions modulo ideals, on the cusp of pcf theory.

\begin{proposition}
\label{prop2}
Let $I$ be an ideal on the cardinal $\kappa$. The following two
statements are equivalent for any regular cardinal $\theta<\kappa$:
\begin{enumerate}
\sk
\item $I$ is weakly $\theta$-saturated and $\theta$-indecomposable.

\sk

\item Suppose $\langle S_\alpha:\alpha<\kappa\rangle$ is a sequence of
sets of ordinals with $|S_\alpha|<\theta$ for all~$\alpha$. Then any
$\leq$-increasing $\theta$-sequence of functions in
$\prod_{\alpha<\kappa}S_\alpha$ is eventually constant modulo $I$.
\sk
\end{enumerate}
\end{proposition}
\begin{proof}
Both directions make use of condition (3) in Proposition~\ref{prop1}.
Let us suppose $\langle S_\alpha:\alpha<\kappa\rangle$ is as in (2), and
assume by way of contradiction that $\bar{h}=\langle
h_i:i<\theta\rangle$ is a $\leq$-increasing sequence in
$\prod_{\alpha<\kappa}S_\alpha$ that is not eventually constant
modulo $I$. By passing to a subsequence we may assume
\begin{equation*}
A_i:=\{\alpha<\kappa:f_i(\alpha)<f_{i+1}(\alpha)\}\notin I\text{ for
all }i<\theta.
\end{equation*}
Fix $\alpha<\kappa$, and suppose $\alpha\in A_i\cap A_j$ for some
$i<j<\theta$.  Since $\bar{h}$ is $\leq$-increasing, we have
\begin{equation}
h_i(\alpha)<h_{i+1}(\alpha)\leq h_j(\alpha).
\end{equation}
Since $|S_\alpha|<\theta$, we conclude
\begin{equation*}
\left|\{i<\theta: \alpha\in A_i\}\right|\leq |S_\alpha|<\theta,
\end{equation*}
and therefore.
\begin{equation}
\bigcap_{i<\theta}\bigcup_{i\leq j<\theta} A_j =\emptyset.
\end{equation}
Each $A_i$ was assumed to be $I$-positive, so we have
contradicted~(1).

For the other direction, let us assume that~(1) fails. By
Proposition~\ref{prop1}, we can find a sequence $\langle
A_i:i<\theta\rangle$ of $I$-positive subsets of~$\kappa$ such that
\begin{equation*}
\bigcap_{i<\theta}\bigcup_{i\leq j<\theta} A_j = \emptyset,
\end{equation*}
that is,
\begin{equation}
\label{eqn3} \left|\{i<\theta:\alpha\in A_i\}\right|<\theta\text{
for all }\alpha<\kappa.
\end{equation}

We define a sequence $\bar{h}=\langle h_i:i<\theta\rangle$ of
functions in $^\kappa\ord$ by the following recursion:

\begin{flushleft}

\noindent{\sf Case 1: Initial stage}

\smallskip

We define  $h_0$ to be identically $0$.

\bigskip

{\sf Case 2: Successor stages}

\smallskip

Given $h_i$, we define $h_{i+1}$ via the formula
\begin{equation*}
h_{i+1}(\alpha) =
\begin{cases}
h_i(\alpha)+1 &\text{if $\alpha\in A_i$}\\
h_i(\alpha)   &\text{otherwise.}
\end{cases}
\end{equation*}

\bigskip
{\sf Case 3: Limit stages}

\medskip

If $i$ is a limit ordinal, then we define $h_{i}$ by setting
$h_i(\alpha) = \sup\{h_j(\alpha):j<i\}$.

\end{flushleft}

Now let us define $S_\alpha:=\{f_i(\alpha):\alpha<\theta\}$. Since
(\ref{eqn3}) holds, our construction guarantees that
$|S_\alpha|<\theta$ for all $\alpha$.  Clearly $\bar{h}$ is
$\leq$-increasing, and since no $A_i$ is in~$I$, the sequence is
also not eventually constant modulo $I$.
\end{proof}

The proof of the above is easily modified to yield the following slightly strengthened result:
\begin{corollary}
Suppose $I$ is weakly $\theta$-saturated and $\theta$-indecomposable ideal on $\kappa$ for some regular cardinal $\theta$,
and let $\langle S_\alpha:\alpha<\kappa\rangle$ be a sequence of sets of ordinals with $|S_\alpha|<\theta$.
If $\bar{h}=\langle h_\beta:\beta<\tau\rangle$ is a $\leq$-increasing sequence of functions in $\prod_{\alpha<\kappa}S_\alpha$
with $\cf(\tau)\geq\theta$, then $\bar{h}$ is eventually constant modulo~$I$.
\end{corollary}

Another easy characterization obtained by slightly different methods is the following:

\begin{proposition}
The following statements are equivalent for an ideal $I$ on $\kappa$
and regular $\theta<\kappa$:
\begin{enumerate}
\item $I$ is weakly $\theta$-saturated and $\theta$-indecomposable.
\sk
\item Any function $f:\kappa\rightarrow\theta$ is bounded below
$\theta$ almost everywhere, that is, there is an ordinal
$\beta<\theta$ such that $f(\alpha)<\beta$ for almost all
$\alpha<\kappa$.
\end{enumerate}
\end{proposition}
\begin{proof}
Assume $f:\kappa\rightarrow\theta$, and for each $\beta<\theta$, let
us define
\begin{equation*}
A_\beta:=\{\alpha<\kappa: f(\alpha)<\beta\}.
\end{equation*}
The sequence $\langle A_\beta:\beta<\theta\rangle$ is increasing
with union $\kappa$, so by (2) of Proposition~\ref{prop1}, there is a
$\beta<\theta$ such that $\kappa\setminus A_\beta \in I$, and
therefore $f(\alpha)<\beta$ for almost all $\alpha<\kappa$.

We prove the other direction by contrapositive --- if (1) fails,
then by (3) of Proposition~\ref{prop1} we can find a family $\langle
A_i:i<\theta\rangle$ of $I$-positive subsets of $\kappa$ such that
\begin{equation}
\label{eqn4} \bigcap_{i<\theta}\bigcup_{i\leq j<\theta} A_j
=\emptyset.
\end{equation}
Now let us define a function $f$ with domain $\kappa$ by
\begin{equation*}
f(\alpha)=\sup\{\beta<\theta: \alpha\in A_\beta\}.
\end{equation*}
Since (\ref{eqn4}) holds, it follows that $f$ maps $\kappa$ into
$\theta$.  However, $f$ is not bounded below~$\theta$ almost
everywhere, as for each $\beta<\theta$, we know
\begin{equation*}
A_\beta\subs\{\alpha<\kappa: \beta\leq f(\alpha)\}
\end{equation*}
and $A_\beta$ is $I$-positive.
\end{proof}

The proof of the above easily yields something slightly stronger:
\begin{corollary}
\label{obviousmod}
Suppose $I$ is a weakly $\theta$-saturated $\theta$-indecomposable ideal on~$\kappa$ for some regular~$\theta<\kappa$.
If $\delta$ is an ordinal of cofinality~$\theta$, then any function $f:\kappa\rightarrow\delta$ is bounded
below~$\delta$ almost everywhere.
\end{corollary}

We close this section by formulating a theorem summarizing the above results.

\begin{theorem}
\label{thm1}
The following statements are equivalent for an ideal $I$ on $\kappa$
and a regular cardinal $\theta<\kappa$.
\begin{enumerate}
\item $I$ is weakly $\theta$-saturated and $\theta$-indecomposable.

\sk

\item Whenever $\langle B_i:i<\theta\rangle$ is an increasing $\theta$-sequence of subsets of $\kappa$, there is an $i^*$ such that
\begin{equation}
i^*\leq i<\theta\Longrightarrow B_i =_I \bigcup_{j<\theta}B_j.
\end{equation}

\sk

\item Whenever $\langle A_i:i<\theta\rangle$ is a $\theta$-sequence of $I$-positive subsets of~$\kappa$, we have
\begin{equation}
\bigcap_{i<\theta}\bigcup_{i\leq j<\theta} A_j \neq \emptyset.
\end{equation}

\sk

\item If $\langle S_\alpha:\alpha<\kappa\rangle$ is a sequence of sets
of ordinals with $|S_\alpha|<\theta$ for all~$\alpha$, then any
$\leq$-increasing $\theta$-sequence of functions in
$\prod_{\alpha<\kappa}S_\alpha$ is eventually constant modulo $I$.

\sk

\item Any function $f:\kappa\rightarrow\theta$ is bounded below
$\theta$ almost everywhere.
\end{enumerate}
\end{theorem}

\section{Least functions and a coloring theorem}

In this section, we show that ideals~$I$ of the sort considered  in Theorem~\ref{thm1} necessarily entail the existence
of a stationary $S^*\subs\kappa$ for which strong versions of $\refl(S^*)$ (involving simultaneous reflection of
stationary subsets of $S^*$) hold.  Once this is established, the section closes by proving a weak version of our main theorem that does not require the complicated machinery introduced later in the paper.
Once again, we need to recall some terminology from the general theory of ideals.

\begin{definition}
Let $I$ be an ideal on a cardinal $\kappa$.
\begin{enumerate}
\item A function $f:\kappa\rightarrow\kappa$ is bounded modulo $I$ if there is a $\xi<\kappa$
such that $\{\alpha<\kappa:\xi<f(\alpha)\}\in I$.

\sk
\item A function $f:\kappa\rightarrow\kappa$ is a {\em least function modulo $I$}
if $f$ is not bounded modulo~$I$, but
\begin{equation*}
g<_I f\Longrightarrow \text{$g$ is bounded modulo $I$}.
\end{equation*}
\end{enumerate}
\end{definition}

The following proposition shows the relevance of the preceding definition to the properties we studied in the preceding section. Shelah (Claim III.3.2A of \cite{cardarith}) obtains the same conclusion from a weaker hypothesis (his proof is the natural generalization of
work of Kanamori and Ketonen  \cite{aki, ketonen1, ketonen2} to the context of not-necessarily-prime ideals), but we include this proof in the interest of completeness and because our assumptions simplify the argument.

\begin{proposition}[Shelah]
\label{shelahlemma}
Let $I$ be an ideal on $\kappa$, and suppose $I$ is
$\theta$-indecomposable and weakly $\theta$-saturated for some
regular $\theta<\kappa$.  Then there is a function $f^*:\kappa\rightarrow\kappa$ such
that
\begin{enumerate}

\sk

\item for each $\beta<\kappa$ we have $f^*(\alpha)>\beta$ for almost all $\alpha<\kappa$, and

\sk

\item if $g(\alpha)<f^*(\alpha)$ for almost all $\alpha<\kappa$, then there is a $\beta<\kappa$ such that
$g(\alpha)<\beta$ for almost all $\alpha<\kappa$.
\end{enumerate}
In particular, $f^*$ is a least function modulo $I$.
\end{proposition}
\begin{proof}
When referring to statement~(1), we will be a little imprecise and say that ``$f^*$ is an upper bound for the constant functions
modulo~$I$'' and rely on context to make it clear to which constant functions we are referring. We will similarly refer
to the conclusion of~(2) by saying ``$g$ is bounded modulo $I$''.

\begin{lemma}
Under our assumptions, if $f$ is an upper bound for the constant functions modulo~$I$ for which~(2) fails, then
we can find a function $h:\kappa\rightarrow\ord$ such that
\begin{itemize}
\sk
\item $h$ is an upper bound for the constant functions modulo~$I$,
\sk

\item $h(\alpha)\leq f(\alpha)$ for all $\alpha<\kappa$, but
\sk
\item $\{\alpha<\kappa:h(\alpha)<f(\alpha)\}\in I^+$.
\sk
\end{itemize}
\end{lemma}
\begin{proof}
By our assumptions, we can find a function $g:\kappa\rightarrow\ord$
such that $g<_I f$, but $g$ is not bounded modulo~$I$. We can freely modify $g$ on sets in $I$, so we may as well
assume $g(\alpha)\leq f(\alpha)$ for all $\alpha<\kappa$.

For each $\beta<\kappa$, let us define
\begin{equation*}
A_\beta=\{\alpha<\kappa: g(\alpha)\leq\beta\}.
\end{equation*}
Since the sequence $\langle A_\beta:\beta<\kappa\rangle$ is $\subs$-increasing, we can apply Proposition~\ref{ob1}
and conclude that there is a $\beta^*$ such that
\begin{equation}
\label{eqn5}
\beta^*\leq\beta<\kappa\rightarrow A_\beta=_I A_{\beta^*}.
\end{equation}
Note that the complement of $A_{\beta^*}$ is not in $I$ because we assumed $g$ is not bounded modulo~$I$. Furthermore,
by~(\ref{eqn5}) we see that
\begin{equation*}
\{\alpha\in \kappa\setminus A_{\beta^*}: g(\alpha)\leq\beta\}\in I\text{ for all }\beta<\kappa.
\end{equation*}
Thus, if we define
\begin{equation*}
h(\alpha)=
\begin{cases}
f(\alpha) &\text{if $\alpha\in A_\beta^*$, and}\\
g(\alpha) &\text{if $\alpha\in\kappa\setminus A_{\beta^*}$},
\end{cases}
\end{equation*}
we have what we need. (Notice that we only needed the weak $\theta$-saturation for the proof of this lemma.)
\end{proof}

Assuming by way of contradiction that there is no function $f^*$ answering conditions (1) and (2), we use the preceding
lemma to build build a $\lvertneqq_I$-decreasing sequence of functions $\langle f_\xi:\xi<\theta\rangle$, each of which is an upper bound
for the constant functions modulo~$I$, using the following inductive recipe:

Let $f_0:\kappa\rightarrow\kappa$ be the identity function.
Since $I$ contains all bounded subsets of~$\kappa$, it
is clear that $f_0$ is an upper bound for the constant functions modulo~$I$.

Given the function $f_\xi$,
we obtain $f_{\xi+1}$ by applying the preceding claim, and for a limit ordinal $\xi<\theta$, we proceed as follows:

For each $\alpha<\kappa$, let us define
\begin{equation*}
S_\alpha:=\{f_\zeta(\alpha):\zeta<\xi\}\cup\{\kappa\},
\end{equation*}
and note that $|S_\alpha|<\theta$ for each $\alpha<\kappa$.
Given  $\beta<\kappa$, we  can define a function $h_\beta\in\prod_{\alpha<\kappa}S_\alpha$ by
\begin{equation*}
h_\beta(\alpha)=\min(S_\alpha\setminus\beta).
\end{equation*}

One should view $h_\beta$ as a ``projection'' of the function that is constant with value~$\beta$ up into the
product $\prod_{\alpha<\kappa}S_\alpha$. Our choice of $f_0$ guarantees that $h_\beta(\alpha)<\kappa$ for almost
all~$\alpha<\kappa$. Finally, we note that
the sequence $\langle h_\beta:\beta<\kappa\rangle$ is $\leq$-increasing, and so an application of Proposition~\ref{prop2} tells us that there is an ordinal $\beta(\xi)<\kappa$ such that $h_\beta =_I h_{\beta(\xi)}$
whenever $\beta(\xi)\leq\beta<\kappa$.  We now define
\begin{equation*}
f_\xi:= h_{\beta(\xi)}.
\end{equation*}

Notice that $f_\xi$ is an upper bound for the constant functions modulo~$I$, and moreover
 for each $\zeta<\xi$, if $\beta(\xi)\leq f_\zeta(\alpha)$ (something that happens for almost all $\alpha<\kappa$), then it must be the case that
\begin{equation}
\label{eqn6}
f_\xi(\alpha)=h_{\beta(\xi)}(\alpha)=\min(S_\alpha\setminus\beta(\xi))\leq f_\zeta(\alpha).
\end{equation}
Thus $f_\xi\leq_I f_\zeta$ for every $\zeta<\xi$.  Since $\xi$ is a limit ordinal, our actions at successor stages
guarantee that in fact $f_\xi\lvertneqq_I f_\zeta$ whenever $\zeta<\xi$.

The preceding construction generates  a sequence of functions $\langle f_\zeta:\zeta<\theta\rangle$ as well as ordinals $\beta(\xi)$ for each limit ordinal $\xi<\theta$. Let us define
\begin{equation*}
\beta^*:=\sup\{\beta(\xi(\epsilon)):\epsilon<\theta\},
\end{equation*}
and for each limit $\xi<\theta$, we set
\begin{equation}
A_\xi :=\{\alpha<\kappa:\beta^*\leq f_{\xi+1}(\alpha)<f_{\xi}(\alpha)\}.
\end{equation}

First, note that each $A_\xi$ is $I$-positive, as $f_{\xi+1}\lvertneqq_I f_\xi$ and $f_{\xi+1}$ is an upper bound
for the constant functions modulo~$I$.  More importantly, we have the following claim:

\begin{claim}
Suppose $\xi<\xi^*<\theta$ are limit ordinals. If $\alpha\in A_\xi$, then $f_{\xi^*}(\alpha)<f_\xi(\alpha)$.
\end{claim}
\begin{proof}
Suppose $\alpha\in A_\xi$, so $\beta^*\leq f_{\xi+1}(\alpha)<f_\xi(\alpha)$. Since $\xi^*$ is a limit ordinal, we
know $\xi+1<\xi^*$. Furthermore, the definition of $\beta^*$ tells us
\begin{equation*}
\beta(\xi^*)\leq f_{\xi+1}(\alpha).
\end{equation*}
From~(\ref{eqn6}), we conclude
\begin{equation*}
f_{\xi^*}(\alpha)\leq f_{\xi+1}(\alpha).
\end{equation*}
Thus, we have
\begin{equation*}
f_{\xi^*}(\alpha)\leq f_{\xi+1}(\alpha)<f_\xi(\alpha),
\end{equation*}
and the proof is complete.
\end{proof}
The above claim makes it clear that no $\alpha<\kappa$ can belong to infinitely many of the sets $A_\xi$.
In particular, this collection of $I$-positive sets contradicts part~(3) of Theorem~\ref{thm1}, and we are done.
\end{proof}

Our next move is to show that the conjunction of $\theta$-indecomposability and weak $\theta$-saturation
has strong consequences for stationary reflection. The following definition will allow us to state some conclusions more precisely:

\begin{definition}
\label{idealterms}
Let $I$ be an ideal on the cardinal~$\kappa$.
\begin{enumerate}
\item $\comp(I)$ is the least cardinal $\tau$ for which $I$ is $\tau$-complete, that is, for which~$I$ is closed under unions of fewer than $\tau$ sets.
\sk
\item $\wsat(I)$ is the least cardinal $\theta$ for which $I$ is weakly $\theta$-saturated.
\sk
\item $\indec(I)=\{\tau<\kappa: \text{$I$ is $\tau$-indecomposable}\}$.
\sk
\item $S^*(I)=\{\alpha<\kappa: \wsat(I)\leq\cf(\alpha)<\alpha\text{ and }\cf(\alpha)\in\indec(I)\}$
\end{enumerate}
\end{definition}

Consider the set $S^*(I)$ defined above for a moment.  Notice that it is non-empty if and only if there is
a regular~$\theta<\kappa$ for which $I$ is both weakly $\theta$-saturated and $\theta$-indecomposable. Also,
if $S^*(I)$ is non-empty, then it is stationary. The omission of regular cardinals from $S^*(I)$ is only relevant
if $\kappa$ happens to be Mahlo; we define the set this way so that Theorem~\ref{thm2} applies uniformly to
any cardinal.

\begin{theorem}
\label{thm2}
Let $I$ be an ideal on the cardinal~$\kappa$. If there is a regular cardinal $\theta<\kappa$
such that $I$ is weakly $\theta$-saturated and $\theta$-indecomposable, then
\begin{enumerate}
\item $S^*(I)$ is stationary,
\sk
\item there is a least function $f^*:\kappa\rightarrow\kappa$ modulo $I$, and
\sk
\item if $S$ is a  stationary subset of $S^*(I)$, then $S\cap f^*(\alpha)$ is stationary in $f^*(\alpha)$
for almost all $\alpha<\kappa$.
\end{enumerate}
In particular,
\begin{enumerate}
\setcounter{enumi}{3}
\item $S^*(I)$ is a stationary subset of $\kappa$ for which $\refl(<\!\comp(T), S^*(I))$ holds.
\end{enumerate}
\end{theorem}
\begin{proof}
We have already remarked that~(1) is a consequence of the given hypotheses, and~(2) is the conclusion of Proposition~\ref{shelahlemma}.
It should also be clear that~(4) follows from~(3), so we will spend our time establishing~(3).

Let $S$ be a stationary subset of $S^*(I)$.  Since~(3) remains true if we establish its conclusion for a stationary subset
of~$S$, we can take advantage of the fact that  cofinality function is regressive on $S^*(I)$ and assume that
$S$ is a subset of $S^\kappa_\tau$ for some $\tau$.  Notice that $I$ is weakly $\tau$-saturated and $\tau$-indecomposable
because of the definition of~$S^*(I)$.

It suffices to prove that whenever we are given a sequence $\langle C_\alpha:\alpha<\kappa, \text{ $\alpha$ limit}\}$ with each $C_\alpha$ closed unbounded in~$\alpha$, we have
\begin{equation}
\label{reflectgoal}
 S\cap C_{f^*(\alpha)}\neq\emptyset \text{ for almost all $\alpha<\kappa$.}
\end{equation}

Given $\beta<\kappa$, let us define $f_\beta:\kappa\rightarrow\kappa$ by
\begin{equation*}
f_\beta(\alpha)=
\begin{cases}
\min(C_{f^*(\alpha)}\setminus\beta) &\text{if $\beta<f^*(\alpha)$,}\\
0 &\text{otherwise.}
\end{cases}
\end{equation*}
Bear in mind that $f_\beta(\alpha)<f^*(\alpha)$ for almost all $\alpha<\kappa$, and so our choice of~$f^*$
implies $f_\beta$ is bounded almost everywhere.  Thus, there is a function $g:\kappa\rightarrow\kappa$
with the property that for any $\beta<\kappa$,
\begin{equation}
\min(C_{f^*(\alpha)}\setminus\beta)<g(\beta)\text{ for almost all }\alpha<\kappa.
\end{equation}

Let $E$ be the closed unbounded subset of~$\kappa$ consisting of those ordinals closed under the function~$g$, and
fix $\delta\in E\cap S$.  Since $\delta$ is closed under~$g$, we know
\begin{equation}
\beta<\delta\Longrightarrow\beta\leq\min(C_{f^*(\alpha)}\setminus\beta)<\delta\text{ for almost all }\alpha<\kappa.
\end{equation}

Now we define another function $h$ by
\begin{equation}
h(\alpha)=
\begin{cases}
\sup(C_{f^*(\alpha)}\cap\delta) &\text{if $\delta\notin C_{f^*(\alpha)}$,}\\
0 &\text{otherwise.}
\end{cases}
\end{equation}
The function $h$ maps $\kappa$ to $\delta$. Since $\cf(\delta)=\tau$ and $I$ is $\tau$-indecomposable, we conclude
from Corollary~\ref{obviousmod} that there is a $\beta<\delta$ for which
\begin{equation}
h(\alpha)<\beta\text{ for almost all }\alpha<\kappa.
\end{equation}
Thus, almost all $\alpha<\kappa$ satisfying the following:
\begin{itemize}
\item $\beta\leq\min(C_{f^*(\alpha)}\setminus\beta)<\delta$, and
\sk
\item $h(\alpha)<\beta$.
\end{itemize}

We finish the proof of~(\ref{reflectgoal}) by establishing that $\delta\in C_{f^*(\alpha)}$ for all~$\alpha<\kappa$
which satisfy both of these statements.  Given such an $\alpha$, assume by way of contradiction that $\delta\notin C_{f^*(\alpha)}$.
On one hand we must have
\begin{equation}
\label{onehand}
h(\alpha)=\sup(C_{f^*(\alpha)}\cap\delta)<\beta,
\end{equation}
while on the other, we have
\begin{equation}
\beta\leq\min(C_{f^*(\alpha)}\setminus\beta)<\delta.
\end{equation}
This latter equation implies
\begin{equation}
\label{theother}
\beta\leq\sup(C_{f^*(\alpha)}\cap\delta),
\end{equation}
and clearly~(\ref{onehand}) and~(\ref{theother}) contradict each other.
\end{proof}

After we noticed the above theorem, we discovered that Shelah uses essentially the same argument in a different context --- it appears tucked into  ``Proof of 3.3 in Case~$\beta$, Subcase (a): \underline{Second Proof}'' on page~149 of~\cite{cardarith} in a result having to do with weakly inaccessible cardinals. Theorem~\ref{thm2} still appears to be new --- one can view it as ``what Shelah's argument really shows".

To see the power of the preceding result, we include the following theorem. This result is essentially a special case of our main theorem. It arises from combining the preceding theorem with one of the main results from~\cite{819}, and its proof serves as a prototype for the argument we employ in Section~10.

\begin{theorem}
\label{wimpythm}
If $\mu$ is a singular cardinal and $\mu^+\rightarrow[\mu^+]^2_{\mu^+}$ holds, then there is a regular $\theta<\mu$ for
which $\refl(<\cf(\mu),S^{\mu^+}_{\geq\theta})$ is true.
\end{theorem}
\begin{proof}
Corollary~5.2 of \cite{819} (which itself relies on results in~\cite{nsbpr}) tells us that under our assumption, there is an ideal $I$ on $\mu^+$ such that
\begin{itemize}
\sk
\item $I$ is $\cf(\mu)$-complete,
\sk
\item $I$ is $\tau$-indecomposable for all regular $\tau$ with $\cf(\mu)<\tau<\mu$, and
\sk
\item $I$ is weakly $\theta$-saturated for some $\theta<\mu$.
\sk
\end{itemize}
Since $\mu^+$ is certainly not a Mahlo cardinal, it follows from the above that
$S^*(I)$ is equal to $S^{\mu^+}_\theta$ modulo the non-stationary ideal, and therefore
 the conclusion follows from Theorem~\ref{thm2}.
\end{proof}

We remark that in the case where $\cf(\mu)$ is uncountable, we can use results from~\cite{nsbpr} to get
 the same conclusion from the failure of $\pr_1(\mu^+,\mu^+,\mu^+,\cf(\mu))$; we will have more to say about this later.

\section{Off-center club guessing}

Our focus now shifts away from the general theory of ideals to questions involving club-guessing.  The proof of Theorem~\ref{wimpythm} turned on properties of an ideal~$I$ whose existence will seem somewhat mysterious to those not familiar with~\cite{nsbpr} and~\cite{819}.
The ideal referenced in the proof of Theorem~\ref{wimpythm} is related to club-guessing, and this explains why we turn
our attention to this matter.

In general, a prototypical club-guessing theorem provides one with a stationary subset $S$ of some cardinal $\kappa$,
 and a sequence $\langle C_\delta:\delta\in S\rangle$ (called an {\em $S$-club sequence})
such that
\begin{itemize}
\item $C_\delta$ is closed and unbounded in $\delta$ for each $\delta\in S$, and
\sk
\item for every closed unbounded $E\subs\kappa$, there are ``many'' $\delta\in S$ for which $E\cap C_\delta$ is ``large''.
\sk
\end{itemize}
One can require various conditions on the sets $C_\delta$, as well as varying the specific meaning of ``many'' and ``large''.

The most well-known club-guessing theorems give us (in certain circumstances) an $S$-club sequence
 $\langle C_\delta:\delta\in S\rangle$  for which $\otp(C_\delta)=\cf(\delta)$, and such that
 for each closed unbounded $E\subs\kappa$, there are stationarily many $\delta\in S$ for which $C_\delta\subs E$.  It
 should be clear that such theorems fit our prototype.

The techniques of \cite{cardarith, 535, nsbpr, 819}, however, require a special sort of club-guessing that is,
in a sense, both weaker and stronger than the standard results. We will give a rough explanation of this, but
we need to enlarge our vocabulary a little bit first:

\begin{definition}
Suppose $C$ is a closed unbounded subset of an ordinal $\delta$. We define
\begin{itemize}
\item $\acc(C)=\{\alpha<\delta: \alpha=\sup(\alpha\cap A)\}\subs C$, and
\sk
\item $\nacc(C)=C\setminus\acc(C)$.
\sk
\end{itemize}
Here ``$\acc$'' stands for ``accumulation points'' and ``$\nacc$'' for non-accumulation points. Finally,
\begin{itemize}
\sk
\item if $\alpha\in\nacc(C)$, then we define $\Gap(\alpha, C)$, {\em the gap in $C$ determined by $\alpha$}, by
\begin{equation}
\Gap(\alpha, C) = (\sup (C\cap\alpha), \alpha).
\end{equation}
\end{itemize}

\end{definition}

Returning now to our discussion, the earlier papers needed results stating that if $\mu$ is singular and $S\subs S^{\mu^+}_{\cf(\mu)}$, then there is an $S$-club sequence $\langle C_\delta:\delta\in S\rangle$ such that for every closed unbounded $E\subs\mu^+$, there are stationarily many $\delta\in S$ such that
\begin{equation*}
(\forall\tau<\mu)\left[\{\alpha\in\nacc(C_\delta)\cap E:\cf(\alpha)>\tau\}\text{ is unbounded in }\delta\right].
\end{equation*}
This type of guessing is weaker than the standard sort in that we aren't requiring $C_\delta$ to be a subset $E$,
but it is also stronger in that we do demand that $C_\delta\cap E$ contains lots of ordinals of large cofinality.
Without going into specifics, if the cofinality of $\mu$ is uncountable, then we can prove the existence of extremely
nice club-guessing sequences with the properties we want (these are the $S$-good pairs of~\cite{nsbpr}) --- sequences
so nice that they can be used to generate colorings of $[\mu^+]^2$.  The club-guessing
result we used in that paper (due originally to Shelah, but a proof can be found in~\cite{819}) does not seem
 to generalize to the case where $\mu$ has countable cofinality. The paper~\cite{819} obtains a weaker result
  for that case, and the conclusions drawn there are correspondingly weaker --- this is why
   Theorem~\ref{wimpythm} refers to square-brackets relations instead
    of the stronger $\pr_1(\mu^+,\mu^+,\mu^+,\cf(\mu))$, and it also explains the existence of the current work.

One of the main goals of this paper is to remedy the situation by moving the club-guessing ``off-center'' --- instead of focusing on subsets of $S^{\mu^+}_{\cf(\mu)}$, we look at stationary subsets consisting of ordinals of larger cofinality.   This is done at a price, for the sets $C_\delta$ we construct are necessarily much more complex. Despite this added complexity,
we are forced to keep tight control over their structure because we still want to be able to connect these
club-guessing sequences with the existence of complicated colorings.

This section presents the club-guessing result alluded to in the preceding paragraph. The theorem is established by
modifying some club-guessing arguments from~\cite{819} to this new context. The main difficulty in this generalization
has to do with the proliferation of parameters, so we start with a list of our main assumptions and notation:

\begin{itemize}
\item $\lambda=\mu^+$ for $\mu$ a singular cardinal
\sk
\item $\kappa=\cf(\mu)$
\sk
\item $\kappa<\sigma=\cf(\sigma)<\mu$
\sk
\item $S$ is a stationary subset of $S^\lambda_\sigma$
\sk
\item $\langle \mu_i:i<\kappa\rangle$ is a continuous increasing sequence of cardinals cofinal in $\mu$,
\sk
\item $\langle c_\delta:\delta\in S\rangle$ is a family of functions such that
\begin{itemize}
\sk
\item $c_\delta$ is an increasing and continuous function from $\sigma$ onto a cofinal subset of $\delta$, and
    \sk
\item for every closed unbounded $E\subs\lambda$, there are stationarily many $\delta\in S$ for which $\ran(c_\delta)\subs E$.
\end{itemize}
\item for $\alpha<\sigma$, $I^\delta_\alpha$ denotes the half-open interval $(c_\delta(\alpha), c_\delta(\alpha+1)]$.
\end{itemize}

Note that if we were to define $C_\delta$ to be $\ran(c_\delta)$ for $\delta\in S$, then we would end up with a
standard sort of club-guessing sequence --- for every closed unbounded $E\subs\lambda$, there would be stationarily
many $\delta\in S$ for which $C_\delta\subs E$.

The next definition captures some standard ideas from proofs of club-guessing; in some cases we have chosen to use more descriptive
names (mostly due to Kojman~\cite{abc}) for these operations than the terminology prevalent in~\cite{cardarith}.

\begin{definition}
Suppose $C$ and $E$ are sets of ordinals with $E\cap \sup(C)$ closed in $\sup(C)$.
We define
\begin{equation}
\drop(C, E) = \{\sup(\alpha\cap E):\alpha\in C\setminus\min(E)+1\}.
\end{equation}
Furthermore, if $C$ and $E$ are both subsets of some cardinal $\lambda$ and $\langle e_\alpha:\alpha<\lambda\rangle$
is a $C$-sequence, then for each $\alpha\in\nacc(C)\cap\acc(E)$, we define
\begin{equation}
\fil(\alpha, C, E) = \drop(e_\alpha, E)\cap\Gap(\alpha, C).
\end{equation}
\end{definition}

The names help one to visualize what the operations do --- $\drop(C, E)$ is the result of ``dropping'' $C$ into the set $E$, while ``$\fil$'' gives us a reasonably canonical way of turning non-accumulation points into accumulation points.
Given these two operations, we are now in a position to prove the following club-guessing theorem:

\begin{theorem}
\label{offcenter}
There is an $S$--club system $\bar{C}=\langle C_\delta:\delta\in S\rangle$ such that
\begin{enumerate}
\item $\ran(c_\delta)\subs C_\delta$
\sk
\item for each $\epsilon<\sigma$ and $i<\kappa$, $C_\delta\cap I^\delta_{\epsilon\cdot\kappa+i}$ has cardinality $\leq\mu_i^+$

\sk

\item if $\alpha\in\nacc(C_\delta)\cap I^\delta_{\epsilon\cdot\kappa+i}$, then $\cf(\alpha)>\mu_i^+$

\sk
\item for every club $E\subs\lambda$, for stationarily many $\delta\in S$, for every $\epsilon<\sigma$ and $i<\kappa$, $E\cap\nacc(C_\delta)\cap I^\delta_{\epsilon\cdot\kappa+i}$ is non-empty.
\end{enumerate}
\end{theorem}

\begin{proof}
We start by simplifying our goal somewhat, by noting that it suffices to produce $\bar{C}$ satisfying~(1),~(2), and
the following modified version of~(4):
\begin{enumerate}
\setcounter{enumi}{3}
\item[$(4)'$] for every club $E\subs\lambda$, for stationarily many $\delta\in S$, for every $\epsilon<\sigma$ and $i<\kappa$, $E\cap\nacc(C_\delta)\cap I^\delta_{\epsilon\cdot\kappa+i}$ contains an ordinal of cofinality $>\mu_i^+$.
\end{enumerate}
Why does this suffice?  Given such an $S$-club system, we  simply throw away those members of $\nacc(C_\delta)$
whose cofinalities are too small to obtain something satisfying~(3), and note that the club-guessing properties
we need are not harmed by this pruning.

 Moving on to the proof, let us assume by way of contradiction that there is no  $S$-club system $\bar{C}$
 satisfying (1), (2), and $(4)'$.  Our aim is to exploit this assumption in order to construct a certain sequence
  $\langle\bar{C}^\zeta:\zeta<\sigma^+\rangle$ of $S$-club systems which will then be used to produce a countable decreasing sequence of ordinals.

Let us agree to say that an $S$-club system {\em satisfies the structural requirements of Theorem~\ref{offcenter}} if conditions (1) and (2) of the conclusion of the theorem hold. We will define objects $E_\zeta$ and $\bar{C}^\zeta=\langle C^\zeta_\delta:\delta\in S\rangle$ by induction on $\zeta<\sigma^+$.  Each $E_\zeta$ will be closed and unbounded
in $\lambda$, while each $\bar{C}^\zeta$ will be an $S$-club system satisfying the structural requirements of Theorem~\ref{offcenter}.  Our convention is that ``stage $\zeta$'' in our construction refers to the process of building $E_{\zeta+1}$ and $\bar{C}^{\zeta+1}$ from $E_\zeta$ and $\bar{C}^\zeta$. Our initial set-up is to take $E_0=\lambda$ and $C^0_\delta=\ran(c_\delta)$ for each $\delta\in S$.

\medskip
\noindent{\sf Stage $\zeta$: Defining $E_{\zeta+1}$ and $\bar{C}^{\zeta+1}$}
\medskip

We assume that our construction furnishes us with an $S$-club system $\bar{C}^\zeta$ satisfying the structural requirements of Theorem~\ref{offcenter}.  Our assumption is that the theorem fails, and so there are closed unbounded subsets $E^0_\zeta$ and $E^1_\zeta$ of $\lambda$
such that for each $\delta\in E^0_\zeta$, we can find $\epsilon<\sigma$ and $i<\kappa$ such that
\begin{equation*}
\alpha\in E^1_\zeta\cap\nacc(C^\zeta_\delta)\cap I^\delta_{\epsilon\cdot\kappa+i}\Longrightarrow \cf(\alpha)\leq \mu^+_i.
\end{equation*}
We define
\begin{equation*}
E_{\zeta+1}:=\acc(E_\zeta\cap E^0_\zeta\cap E^1_\zeta).
\end{equation*}

The definition of $\bar{C}^{\zeta+1}$ will take a bit more effort.  Let us agree to call an ordinal $\delta\in S$ {\em active at stage $\zeta$} if $C^0_\delta\subs\acc(E_{\zeta+1})$.  Our choice of $\langle c_\delta:\delta\in S\rangle$ ensures that at any stage, the set of active $\delta$ is a stationary subset of $S$.  If $\delta\in S$ is inactive at stage $\zeta$, then we simply define $C^{\zeta+1}_\delta$ to be $C^\zeta_\delta$ and do nothing.

On the other hand, if $\delta$ {\em is} active at stage $\zeta$, then $\delta$ must be in $E^0_\zeta$ and therefore we can find a least $a(\delta,\zeta)<\sigma$ such that
\begin{equation*}
I(\delta,\zeta):= I^\delta_{a(\delta,\zeta))},
\end{equation*}
we have
\begin{equation}
\label{4.1}
\alpha\in E^1_\zeta\cap\nacc(C^\zeta_\delta)\cap I(\delta,\zeta)\Longrightarrow \cf(\alpha)\leq \mu_{i(\delta,\zeta)}^+.
\end{equation}
Our construction of $C^{\zeta+1}_\delta$ will modify $C^\zeta_\delta$ only on the interval $I(\delta,\zeta)$ --- everything else will be left untouched.

The ordinal $a(\delta,\zeta)$ can be written in the form
\begin{equation}
a(\delta,\zeta)=\epsilon(\delta,\zeta)\cdot\kappa+i(\delta,\zeta)
\end{equation}
for some unique $\epsilon(\delta,\zeta)<\sigma$ and $i(\delta,\zeta)<\kappa$. These two ordinals will also play a role in our construction.

Our next move is to define
\begin{equation*}
D^\zeta_\delta:=\drop(C^\zeta_\delta\cap I(\delta,\zeta), E_{\zeta+1}\cap I(\delta,\zeta)).
\end{equation*}
We note the following facts about $D^\zeta_\delta:$
\begin{itemize}
\item $c_\delta(a(\delta,\zeta)+1)$ --- the top of the interval $I(\delta,\zeta)$ --- is an element of $D^\zeta_\delta$ because $C^0_\delta\subs\acc(E_{\zeta+1})$,
    \sk

\item $D^\zeta_\delta$ is a closed subset of $E_{\zeta+1}\cap I(\delta,\zeta)$,
\sk
\item $|D^\zeta_\delta|\leq |C^\zeta_\delta\cap I(\delta,\zeta)|\leq\mu^+_{i(\delta,\zeta)}$ (as $\bar{C}^\zeta$ satisfies the structural requirements of Theorem~\ref{offcenter}), and
    \sk
\item if $C^\zeta_\delta\cap I(\delta,\zeta)$ is unbounded in $c_\delta(a(\delta,\zeta)+1)$, then so is $D^\zeta_\delta$.
\sk
\end{itemize}
One should picture $D^\zeta_\delta$ as arising after ``shifting'' $C^\zeta_\delta\cap I(\delta,\zeta)$ so that it lies inside of $E_{\zeta+1}$.  Our construction will ensure that $D^\zeta_\delta\subs C^{\zeta+1}_\delta$, but we need to do more work first.

Let us say that an element $\alpha$ of $D^\zeta_\delta$ {\em needs attention} if
\begin{equation}
\alpha\in\acc(E_{\zeta+1})\cap \nacc(D^\zeta_\delta),
\end{equation}
and
\begin{equation}
\cf(\alpha)\leq \mu_{i(\delta,\zeta)}^+.
\end{equation}
If $\alpha$ needs attention, then  $\Fill(\alpha, C^\zeta_\delta\cap I(\delta,\zeta), E_{\zeta+1}\cap I(\delta,\zeta))$ provides us with a closed unbounded subset of $\alpha$ lying in the interval $\Gap(\alpha,D^\zeta_\alpha)$. Notice as well that this closed unbounded subset of $\alpha$ is of cardinality $\cf(\alpha)\leq \mu_{i(\delta,\zeta)}^+$, and so the set
\begin{equation*}
A^\zeta_\delta:=D^\zeta_\delta\cup\{\Fill(\alpha, C^\zeta_\delta\cap I(\delta, \zeta),E_{\zeta+1}\cap I(\delta,\zeta)):\alpha\text{ needs attention}\}
\end{equation*}
satisfies
\begin{equation}
|A^\zeta_\delta|\leq |D^\zeta_\delta|\cdot\mu_{i(\delta,\zeta)}^+\leq |C^\zeta_\delta\cap I(\delta,\zeta)|\cdot\mu_{i(\delta,\zeta)}^+ =\mu_{i(\delta,\zeta)}^+.
\end{equation}
Since the needed instances of ``$\Fill$'' are always closed sets lying in a ``gap'' of $D_\delta^\zeta$, the set $A^\zeta_\delta$ is also closed in $c_\delta(a(\delta,\zeta)+1)$.

Moreover, $A^\zeta_\delta$ is also unbounded in $c_\delta(a(\delta,\zeta)+1)$, for either $D^\zeta_\delta$ is already unbounded, or it is the case that $c_\delta(a(\delta,\zeta)+1)$ itself needs attention.

We now define $C^{\zeta+1}_\delta$ in piecewise fashion:
\begin{equation*}
C^{\zeta+1}_\delta\setminus I(\delta,\zeta)= C^\zeta_\delta\setminus I(\delta,\zeta),
\end{equation*}
and
\begin{equation*}
C^{\zeta+1}_\delta\cap I(\delta,\zeta)= A^\zeta_\delta.
\end{equation*}
The $S$-club system $\bar{C}^{\zeta+1}_\delta$ satisfies the structural requirements of our theorem, and so the construction can continue.

We still need to describe how to obtain $\bar{C}^\zeta$ and $E_\zeta$ when $\zeta$ is a limit.  Our construction defines
\begin{equation*}
E_\zeta=\bigcap_{\xi<\zeta}C_\xi,
\end{equation*}
and for each $\delta\in S$, we define $C^\zeta_\delta$ to be the closure in $\delta$ of
\begin{equation*}
\{\alpha<\delta:\alpha\in C^\xi_\delta\text{ for all sufficiently large }\xi<\zeta\}.
\end{equation*}
Note that $C^\zeta_\delta$ is closed in $\delta$, and it is unbounded in $\delta$ as it contains $C^0_\delta$.  Elementary cardinal arithmetic implies that $\bar{C}_\zeta$ satisfies the structural requirements of Theorem~\ref{offcenter}, and so our construction can continue.

To this point in the proof, we have used the failure of Theorem~\ref{offcenter} to produce a sequence $\langle \bar{C}^\zeta:\zeta<\sigma^+\rangle$ of $S$-club systems; our task is to show that this leads to a contradiction.

Let us define
\begin{equation*}
E^*:=\bigcap_{\zeta<\sigma^+} E_\zeta.
\end{equation*}
Since $E^*$ is club in $\lambda$, we can find a $\delta\in S$ for which
\begin{equation}
\label{eqn9}
C^0_\delta\subs\{\alpha<\lambda:\mu \text{ divides }\otp(E^*\cap\alpha)\},
\end{equation}
and this guarantees that
\begin{equation}
\label{eqn10}
\alpha\in\nacc(C^0_\delta)\Longrightarrow\left|E^*\cap\Gap(\alpha, C^0_\delta)\right|=\mu.
\end{equation}

We know that this $\delta$ is active at each stage $\zeta<\sigma^+$ because of (\ref{eqn9}), and therefore we can find $\epsilon^*<\sigma$ and $i^*<\kappa$ such that $\epsilon(\delta,\zeta)=\epsilon^*$ and $i(\delta,\zeta)=i^*$ for unboundedly many $\zeta<\sigma^+$. Letting
\begin{equation*}
I^*:= (c_\delta(\epsilon^*\cdot\kappa+i^*),c_\delta(\epsilon^*\cdot\kappa+i^*+1)],
 \end{equation*}
this means $I(\delta,\zeta)=I^*$ for unboundedly many ordinals $\zeta<\sigma^+$.  Let $\langle \zeta_n:n<\omega\rangle$ be the increasing enumeration of the first $\omega$ such ordinals, and let $\zeta^*=\sup\{\zeta_n:n\in\omega\}$.

Our construction ensures $|I^*\cap C^\zeta_\delta|\leq \mu_{i^*}^+$ for all $\zeta<\sigma^+$, so an appeal to (\ref{eqn10}) allows us to choose
\begin{equation}
\label{eqn11}
\beta^*\in E^*\cap I^*\setminus\bigcup_{\zeta<\sigma^+} C^\zeta_\delta.
\end{equation}
Finally, define
\begin{equation*}
\beta_n:=\min(C^{\zeta_n}_\delta\setminus\beta^*).
\end{equation*}
Note that our choice of $\beta^*$ ensures that $\beta^*<\beta_n$ for all $n$.

\begin{claim}
For each $n$, we have $\beta_{n+1}<\beta_n$.
\end{claim}
\begin{proof}
Given $n$, we note that
\begin{equation*}
\beta_{n+1}=\min(C^{\zeta_n+1}_\delta)
\end{equation*}
(notice the shift in the position of ``$+1$'' here), as the definition of $\zeta_{n+1}$ implies $C^\xi_\delta\cap I^*=C^{\zeta_n+1}_\delta\cap I^*$ whenever $\zeta_n+1\leq\xi\leq\zeta_{n+1}$.  We now split up into two cases:

\medskip

\noindent{\bf\sf Case 1:} $\beta_n\notin \acc(E_{\zeta_n+1})$.

\medskip

Since $\beta^*\in E^*\subs E_{\zeta_n+1}$, it follows that
\begin{equation*}
\beta^*\leq\sup(\beta_n\cap E_{\zeta_n+1})<\beta_n.
\end{equation*}
Now $\beta^*$ is not in $C^{\zeta_n+1}_\delta$ while
\begin{equation*}
\sup(\beta_n\cap E_{\zeta_n+1})\in D^{\zeta_n+1}_\delta\subs C^{\zeta_n+1}_\delta
\end{equation*}
so $\beta^*<\beta_{n+1}<\beta_n$ as claimed.

\medskip

\noindent{\bf\sf Case 2:} $\beta_n\in\acc(E_{\zeta_n+1})$.

\medskip

Both $\delta$ and $\beta_n$ are in $E_{\zeta_n+1}$, so in particular we know $\delta\in E^0_{\zeta_n}$ and $\beta_n\in E^1_{\zeta_n}$. In addition, $\beta_n$ must lie in $\nacc(C^{\zeta_n}_\delta)$ because $\beta^*<\beta_n$. By (\ref{4.1}), we conclude that $\cf(\beta_n)\leq\mu_{i(\delta,\zeta_n)}^+$.

We have assumed in our case hypothesis that $\beta_n=\sup(E_{\zeta_n+1}\cap\beta_n)$ and this guarantees that $\beta_n$ is an element of $D^{\zeta_n}_\delta$.  Since $\beta_n$ is also in $\nacc(C^{\zeta_n}_\delta)$, the set $D^{\zeta_n}_\delta$ cannot
pick up any new elements between $\beta^*$ and $\beta_n$ and hence
\begin{equation}
\beta_n =\min(D^{\zeta_n}_\delta\setminus\beta^*)>\beta^*.
\end{equation}
Thus, $\beta_n$ is in $\nacc(D^{\zeta_n})$ and we see that $\beta_n$ needs attention during the construction of $C^{\zeta_n+1}_\delta$.

In this case, our construction makes sure that $C^{\zeta_n+1}_\delta$ contains a closed unbounded subset of $\beta_n$.
In particular, $C^{\zeta_n+1}_\delta\cap (\beta^*,\beta_n)\neq\emptyset$,  and therefore
\begin{equation*}
\beta^*<\beta_{n+1}=\min(C^{\zeta_{n+1}}_\delta\setminus\beta^*)=\min(C^{\zeta_n+1}\setminus\beta^*)<\beta_n,
\end{equation*}
as required.
\end{proof}

In summary, if the conclusion of our theorem fails, then our construction generates an infinite decreasing sequence of ordinals. This is absurd, and so the theorem is established.
\end{proof}

\section{Organizational interlude}

Our goal in this section is to lay a good foundation before proceeding to the proofs of our main theorems. A good deal of the difficulty
in these proofs lies in the fact that we require so many different objects to push the argument through; we make an attempt at organizing
notation and providing a clear picture of our assumptions before proceeding.  In addition, we focus only on the case
 where~$\cf(\mu)=\aleph_0$, as uncountable cofinalities can be handled by earlier work. We first fix our names for the
 various cardinals that are important for our theorems:

\begin{itemize}
\item $\lambda=\mu^+$ for $\mu$ singular of cofinality $\aleph_0$.
\sk
\item $\aleph_0<\sigma=\cf(\sigma)<\mu$,
\sk
\item $S$ is a stationary subset of $\{\delta<\lambda:\cf(\delta)=\sigma\}$
\sk
\item $\langle \mu_i:i<\omega\rangle$ is a strictly increasing sequence of regular cardinals such that
\begin{itemize}
\item $\langle \mu_i:i<\omega\rangle$ is cofinal in $\mu$, and
\sk
\item $\sigma<\mu_0$.
\end{itemize}
\end{itemize}

Now the club-guessing result of the previous section gives us objects
 $\langle c^0_\delta:\delta\in S\rangle$,  $\langle C_\delta:\delta\in S\rangle$, and $\{I(\delta, \epsilon, m):\delta\in S, \epsilon<\sigma, m<\omega\}$ satisfying the following:

\begin{itemize}
\item $c^0_\delta$ is the increasing enumeration of a closed unbounded subset of~$\delta$ of order-type~$\sigma$
\sk
\item if $E$ is a closed unbounded subset of $\lambda$, then the set of $\delta\in S$ for which $\ran(c^0_\delta)\subs E$ is stationary
\sk
\item $C_\delta$ is a closed unbounded subset of $\delta$ with $\ran(c^0_\delta)\subs C_\delta$
\sk
\item $I(\delta,\epsilon, m)$ denotes the half-open interval $(c^0_\delta(\omega\cdot\epsilon+m), c^0_\delta(\omega\cdot\epsilon+m+1)]$
\sk
\item $|C_\delta\cap I(\delta,\epsilon, m)|\leq \mu_m^+$,
\sk
\item $\alpha\in \nacc(C_\delta)\cap I(,\delta,\epsilon, m)\Longrightarrow \cf(\alpha)>\mu_m^+$, and
\sk
\item if $E$ is a closed unbounded subset of $\lambda$, then for stationarily many $\delta\in S$ it is the case that for each
 $\epsilon<\sigma$ and $m<\omega$, the set $E\cap\nacc(C_\delta)\cap I(\delta,\epsilon, m)$ is non-empty.
\end{itemize}
If $\delta\in S$ is fixed and clear from context (as is usually the case), then we will  write $I(\epsilon, m)$ instead of $I(\delta,\epsilon, m)$.

It is crucial for the reader to have a good picture of the structure of the objects described above. Our notation
is intended to describe something fairly simple: given $\delta\in S$, we use $c_\delta^0$ (or rather, the range of $c_\delta^0$)
essentially to divide the ordinals less than $\delta$ into $\sigma$ blocks, each of which is further divided into $\omega$ pieces.
We say ``essentially'', because there are a few ordinals left out for technical reasons, but the reader will be well-served
by thinking of the interval $I(\delta, \epsilon, m)$ as ``the $m$th piece in block $\epsilon$ built using $C_\delta$''.  The following observation tells the complete story.

\begin{proposition}
Given $\delta\in S$, each ordinal $\alpha<\delta$ satisfies exactly one of the following conditions:
\begin{enumerate}
\sk
\item $\alpha\leq c_\delta^0(0)$,
\sk
\item $\alpha\in\acc(\ran(c^0_\delta))$, or
\sk
\item $\alpha\in I(\epsilon, m)$ for some unique $\epsilon<\sigma$ and $m<\omega$.
\sk
\end{enumerate}
In particular,any element of $\nacc(C_\delta)\setminus\{\min(C_\delta)\}$ lies in $I(\epsilon, m)$ for some unique
$\epsilon<\sigma$ and $m<\omega$.
\end{proposition}

Our next task is to  connect our club-guessing sequence to some ideals introduced by Shelah in~\cite{cardarith}.
We have already discussed club-guessing in general terms at the start of the previous section; what follows is a
deeper discussion the particular case of interest to us.

\begin{definition}
\label{ideltadef}
Given $\delta\in S$, let $I_\delta$ be the ideal on $C_\delta$ generated by the sets of the form
\begin{equation}
\label{ideltaset}
\{\gamma\in C_\delta: \gamma\in\acc(C_\delta)\text{ or }\cf(\gamma)<\alpha\text{ or }\gamma<\beta\}
\end{equation}
for $\alpha<\mu$ and $\beta<\delta$.
\end{definition}

Shelah uses the notation $J^{b[\mu]}_{C_\delta}$ for the ideal described above; we have gone with something a little simpler. These ideals $I_\delta$ are natural in our given context --- for example, the following two propositions show us that the ideals are compatible with the interval structure we have imposed on the objects $C_\delta$

\begin{proposition}
\label{almostallprop}
Suppose $\delta\in S$.  If $\epsilon^*<\sigma$ and $m^*<\omega$ are fixed, then $I_\delta$-almost members of $C_\delta$ are in $\nacc(C_\delta)\cap I(\epsilon, m)$ for some $\epsilon\geq\epsilon^*$ and $m\geq m^*$.
\end{proposition}
\begin{proof}
This is clear from the definition of $I_\delta$.
\end{proof}

\begin{proposition}
\label{propguess}
If $E$ is a closed unbounded subset of $\lambda$, then there are stationarily many $\delta$ for which $E\cap C_\delta\notin I_\delta$.
\end{proposition}
\begin{proof}
Our assumptions give us  a stationary set of $\delta$ such that $E\cap \nacc(C_\delta)\cap I(\epsilon, m)$ is non-empty for each
$\epsilon<\sigma$ and $m<\omega$, so we are done if we verify that $E\cap C_\delta\notin I_\delta$ for each such $\delta$.
This is easily done --- given $\alpha<\mu$ and $\beta<\delta$, we need to show that $E\cap C_\delta$ contains a member
of $\nacc(C_\delta)$ larger than $\alpha$ and whose cofinality is greater than $\beta$.  Choose $\epsilon<\sigma$ so large
that $\alpha<c_\delta^0(\omega\cdot\epsilon)$ and choose $m$ so large that $\beta<\mu_m$. We know $E\cap\nacc(C_\delta)\cap I(\epsilon, m)$
is non-empty, and any ordinal in this intersection is necessarily a member of $E\cap\nacc(C_\delta)$ larger
than $\alpha$ and of cofinality greater than $\beta$.
\end{proof}

Given the preceding proposition, we can now bring in a certain club-guessing ideal which will be woven into all of our most important results. The following definition makes sense in a more general context than that which we are considering (Chapter III of~\cite{cardarith}, for example), but we will deal only with the particular case of interest to us.

\begin{definition}
\label{idpcidef}
Let $\bar{I}=\langle I_\delta:\delta\in S\rangle$ be as in Definition~\ref{ideltadef}.  The ideal $\id_p(\bar{C},\bar{I})$ is defined by $A\in \id_p(\bar{C},\bar{I})$ if there is a closed unbounded $E\subs\lambda$
such that $$A\cap E\cap C_\delta\in I_\delta \text{ for all }\delta\in S\cap E.$$
\end{definition}

The preceding definition is robust under slight modifications. For example, we get the same ideal by requiring
$A\cap E\cap C_\delta$ to be in $I_\delta$ for all but non-stationarily many $\delta\in S\cap E$. Note as well that the restriction to $\delta\in S\cap E$ is not important --- if $\delta\in S\setminus E$, then $E\cap C_\delta$ must be bounded and therefore $E\cap C_\delta\in I_\delta$.  Roughly speaking, a
set $A\subs\lambda$ is in $\id_p(\bar{C},\bar{I})$ if $\bar{C}$ fails to ``guess'' $E\cap A$ in the sense of
Proposition~\ref{propguess}. Note that in light of Proposition~\ref{propguess}, the ideal $\id_p(\bar{C},\bar{I})$ is a proper ideal on $\lambda$ extending the non-stationary ideal.

The following fact is a specific instance of Observation~3.2 in Chapter~III of Shelah --- it is an elementary result establishing indecomposability properties for the club-guessing ideal under consideration. We include the proof because the result is crucial piece of our main theorem.

\begin{proposition}
\label{indecprop}
Under our assumptions, the ideal $\id_p(\bar{C},\bar{I})$ is $\tau$-indecomposable whenever $\tau<\mu$ is an uncountable regular cardinal distinct from $\sigma$.
\end{proposition}
\begin{proof}
 Let $\tau\neq\sigma$ be an uncountable regular cardinal less than~$\mu$. Given $\delta\in S$, we first establish that $I_\delta$ is $\tau$-indecomposable, so let $\langle A_i:i<\tau\rangle$ be an increasing sequence of subsets of $C_\delta$, each
of which is in $I_\delta$. The collection of sets of the form~(\ref{ideltaset}) is closed under finite unions, and since these sets generate $I_\delta$, it follows that each set in the ideal is covered by a single set of that form.  Thus, we may define $\alpha_i$ to be the least $\alpha<\mu$ for which there is a $\beta<\delta$ such that
\begin{equation*}
A_i\subs\{\gamma\in C_\delta:\gamma\in\acc(C_\delta)\text{ or }\cf(\gamma)<\alpha\text{ or }\gamma<\beta\}.
\end{equation*}
Since the sequence $\langle A_i:i<\tau\rangle$ is increasing, it is clear that $\langle \alpha_i:i<\tau\rangle$ is
a non-decreasing sequence of ordinals less than~$\mu$.  Since $\cf(\mu)=\aleph_0$ and $\tau$ is an uncountable regular cardinal, it follows this sequence is bounded by some ordinal $\alpha^*<\mu$.

Next, we let $\beta_i$ be the least ordinal $\beta<\delta$ for which
\begin{equation*}
A_i\subs\{\gamma\in C_\delta:\gamma\in\acc(C_\delta)\text{ or }\cf(\gamma)<\alpha^*\text{ or }\gamma<\beta\}.
\end{equation*}
Such a $\beta$ always exists because of our choice of $\alpha^*$, and since $\langle A_i:i<\tau\rangle$ is increasing,
the sequence $\langle \beta_i:i<\tau\rangle$ is non-decreasing.  Since $\tau\neq\sigma=\cf(\delta)$, there is an ordinal $\beta^*<\delta$ greater than all $\beta_i$.

Thus, we conclude
\begin{equation}
\bigcup_{i<\tau} A_i \subs\{\gamma\in C_\delta:\gamma\in\acc(C_\delta)\text{ or }\cf(\gamma)<\alpha^*\text{ or }\gamma<\beta^*\},
\end{equation}
and hence this union is in $I_\delta$.

The rest of the argument consists in noting that $\id_p(\bar{C},\bar{I})$ inherits the indecomposability properties
of the ideals $I_\delta$.  To see this, let $\langle A_i:i<\tau\rangle$ be an increasing family of subsets of $\lambda$, each of which is in $\id_p(\bar{C},\bar{I})$.  Set $A$ equal to the union of these sets, and suppose by way of contradiction that $A\notin \id_p(\bar{C},\bar{I})$.

Since $A_i\in \id_p(\bar{C},\bar{I})$, there is a closed unbounded $E_i\subs\lambda$ for which
\begin{equation}
\label{small}
A_i\cap E_i\cap C_\delta\in I_\delta\text{ for all }\delta\in S\cap E.
\end{equation}
If we let $E$ be the intersection of all the sets $E_i$, then it is clear that $E$ is closed and unbounded in $\lambda$. We have assumed that $A\notin\id_p(\bar{C},\bar{I})$, and thus we can fix an ordinal $\delta\in S\cap E$ such that
\begin{equation}
\label{large}
A\cap E\cap C_\delta\notin I_\delta.
\end{equation}
Define
\begin{equation*}
B_i:= A_i\cap E\cap C_\delta.
\end{equation*}
It should be clear that the sequence $\langle B_i:i<\tau\rangle$ is a $\subs$-increasing sequence of subsets of~$C_\delta$, and since $E\subs E_i$, we know $B_i\in I_\delta$ by~(\ref{small}). Now $I_\delta$ is $\tau$-indecomposable, and so
\begin{equation*}
B:=\bigcup_{i<\tau}B_i \in I_\delta.
\end{equation*}
This is a contradiction, though, because we also know
\begin{equation*}
B = \bigcup_{i<\tau}(A_i\cap E\cap C_\delta) = \left(\bigcup_{i<\tau}A_i\right)\cap E\cap C_\delta = A\cap E\cap C_\delta,
\end{equation*}
which is not in $I_\delta$ by~(\ref{large}).  We conclude that $A$ must be in $\id_p(\bar{C},\bar{I})$, as required.
\end{proof}

Now that we have established some basic results about our club-guessing sequence, our concern shifts
to seeing how well some of the techniques used in~\cite{819} can be adapted to this context. We begin
with an observation that the interval structure we have placed on our sets $C_\delta$ gives us a natural way of writing each $C_\delta$ as a countable increasing union of approximations:

\begin{definition}
Given $\delta\in S$ and $m<\omega$, we define
\begin{equation*}
C_\delta[m] = \ran(c^0_\delta)\cup\bigcup\{C_\delta\cap I(\epsilon, i):\epsilon<\sigma\text{ and }i\leq m\}.
\end{equation*}
\end{definition}

In more descriptive language, $C_\delta[m]$ consists of $\ran(c^0_\delta)$ together with the parts of $C_\delta$ lying in the first $m+1$ pieces of each of the $\sigma$ blocks built using $C_\delta$. We note the following easy facts about these objects:

\begin{proposition}
Suppose $\delta\in S$. Then
\begin{enumerate}
\item $C_\delta[m]$ is closed and unbounded in $\delta$,
\sk
\item $|C_\delta[m]|\leq\sigma\cdot\mu_m^+=\mu_m^+$,
\sk
\item $C_\delta[m]\subs C_\delta[m+1]$, and
\sk
\item $C_\delta=\bigcup_{m<\omega} C_\delta[m]$.
\end{enumerate}
\end{proposition}

Our next move is implement ``Shelah's ladder-swallowing trick'' (a crucial ingredient of the proofs in \cite{nsbpr} and \cite{819}) in this new context. We take the following definition from~\cite{819}.

\begin{definition}
 A {\em generalized $C$-sequence} on $\lambda$ is a family $$\langle e_\alpha^m:\alpha<\lambda,m<\omega\rangle$$
such that for each $\alpha<\lambda$ and $m<\omega$,
\begin{itemize}
\item $e^m_\alpha$ is closed unbounded in $\alpha$, and
\medskip
\item $e^m_\alpha\subs e^{m+1}_\alpha$.
\end{itemize}
\end{definition}

We now use our collection of objects $C_\delta[m]$ in order to construct a very special generalized $C$-sequence.

\begin{lemma}
\label{swallowlemma}
There is a generalized $C$-sequence $\langle e^m_\alpha:\alpha<\lambda, m<\omega\rangle$ such that
\begin{enumerate}
\item $|e^m_\alpha|\leq \cf(\alpha)+\mu_m^+$, and
\sk
\item $\delta\in S\cap e^m_\alpha\Longrightarrow C_\delta[m]\subs  e^m_\alpha$
\sk
\end{enumerate}
\end{lemma}
\begin{proof}
Let $e_\alpha$ be closed unbounded in $\alpha$ of order-type~$\cf(\alpha)$.  We define
\begin{eqnarray*}
e^0_\alpha[0]&= &e_\alpha\\
e^0_\alpha[k+1]&=&\text{ closure in $\alpha$ of }e^0_\alpha[k]\cup\bigcup\{C_\delta[0]:\delta\in S\cap e^0_\alpha[k]\}\\
e^0_\alpha & = &\text{ closure in $\alpha$ of }\cup\{e^0_\alpha[k]:k<\omega\}\\
e^{m+1}_\alpha[0] & = &e^m_\alpha\\
e^{m+1}_\alpha[k+1] &= &\text{ closure in $\alpha$ of }e^{m+1}_\alpha[k]\cup \bigcup\{C_\delta[m+1]:\delta\in S\cap e^{m+1}\alpha[k]\}\\
e^m_\alpha &= &\text{ closure in $\alpha$ of }\cup\{e^m_\alpha[k]:k<\omega\}.
\end{eqnarray*}
The estimate for $|e^m_\alpha|$ holds because of the corresponding bounds on the size of each $C_\delta[m]$. As for the other requirement, note that since $\sigma$ is an uncountable regular cardinal,  if $\delta\in S\cap e^m_\alpha$ then there is a $k<\omega$ such that $\delta\in S\cap e^m_\alpha[k]$.
\end{proof}

Both conditions (1) and (2) in the preceding lemma are crucial ingredients in the proof of our main result. The name ``ladder-swallowing trick'' comes from condition~(2), which tells us that the generalized $C$-sequence we have built ``swallows'' the objects $C_\delta[n]$ in a natural way.  These sorts of requirements probably seem quite unmotivated at this point, but we note that it is precisely this ``swallowing'' behavior that allows us to lift some of Todor{\v{c}}evi{\'c}'s minimal walks machinery to successors of singular cardinals without the need for strong assumptions like $\square$ (a technique first pioneered by Shelah in Chapter~III of~\cite{cardarith}, as well as~\cite{535}).

\section{Minimal walks and a preliminary result}

We leave club-guessing behind for a bit, and turn now to Todor{\v{c}}evi{\'c}'s technique of minimal walks as it relates to our specific generalized $C$-sequence.  The following definition presents the
notation we use to keep track of everything.

\begin{definition}
Given ordinals $\alpha<\beta<\lambda$ and $m<\omega$, we define
\begin{equation*}
\beta^m_0(\alpha,\beta)=\beta,
\end{equation*}
while
\begin{equation*}
\beta^m_{i+1}(\alpha,\beta)=
\begin{cases}
\min\left(e^m_{\beta^m_i(\alpha,\beta)}\setminus\alpha\right)  &\text{if $\beta^m_i(\alpha,\beta)>\alpha$, and }\\
\alpha  &\text{otherwise}
\end{cases}
\end{equation*}
Next, we define
\begin{equation*}
\rho^m_2(\alpha,\beta)=\text{ least $i$ for which }\beta^m_i(\alpha,\beta)=\alpha.
\end{equation*}
\end{definition}

For each $m<\omega$ and $\alpha<\beta<\lambda$, the sequence $\langle \beta^m_i:i<\rho^m_2(\alpha,\beta)\rangle$ is the minimal walk from $\beta$ to $\alpha$ using the $C$-system $\bar{e}^m$ in the sense of~\cite{stevobook}. We will abbreviate this somewhat, and call this {\em the $m$-walk from $\beta$ to $\alpha$}. The use of $\rho_2$ to stand for the length of such a walk is due to Todor{\v{c}}evi{\'c}.  Minimal walks in the context of generalized $C$-sequences were introduced in~\cite{819}; our notation is a variant of the notation in that paper.

  Now that we have the definition of $m$-walk at our disposal, the statement of the following theorem makes sense.  This theorem illustrates how the structure of our generalized $C$-system causes lets us establish connection between the $m$-walks defined above, and the club-guessing ideal $\id_p(\bar{C},\bar{I})$ discussed in the preceding section. This theorem will quickly be superseded by results in the next section, but it represents an important step on the way to our main theorem.

 \begin{theorem}
 \label{filterthm}
 Let $\langle t_\alpha:\alpha<\lambda\rangle$ be a pairwise disjoint family of finite subsets of~$\lambda$.  For $\id_p(\bar{C},\bar{I})$-almost all $\beta^*<\lambda$, there are $\alpha<\beta<\lambda$  and $m<\omega$ such that for every $\zeta\in t_\alpha$ and $\xi\in t_\beta$, the $m$-walk from $\xi$ to $\zeta$ passes through $\beta^*$.
 \end{theorem}

The proof of Theorem~\ref{filterthm} will fill the rest of this section. We begin following
lemma and corollary, which are translations into the context of generalized $C$-sequences of some standard facts about minimal walks.  Here we see our assumption that elements of $S$ have uncountable cofinality starts to become relevant.

\begin{lemma}
Suppose $\delta<\beta<\lambda$ and $\cf(\delta)>\aleph_0$. There is an ordinal $\gamma^*(\delta,\beta)<\delta$ such that
\begin{equation}
(\gamma^*(\delta,\beta),\delta)\cap e^m_{\beta^m_i(\delta,\beta)}=\emptyset
\end{equation}
for all $m<\omega$ and $i<\rho^m_2(\delta,\beta)-1$.
\end{lemma}
\begin{proof}
If $m<\omega$ and $i<\rho^m_2(\delta,\beta)-1$, then $\delta\notin e^m_{\beta^m_i(\delta,\beta)}$. Since this latter set is closed, it follows that
\begin{equation*}
\sup(\delta\cap e^m_{\beta^m_i(\delta,\beta)})<\delta,
\end{equation*}
and if we define
\begin{equation*}
\gamma^*(\delta,\beta):=\sup\{\sup(\delta\cap e^m_{\beta^m_i(\delta,\beta)}):m<\omega, i<\rho^m_2(\delta,\beta)-1\},
\end{equation*}
then $\gamma^*(\delta,\beta)<\delta$ as $\delta$ is of uncountable cofinality.
\end{proof}

The following corollary isolates the importance of $\gamma^*(\delta,\beta)$ --- if $\alpha$ is any ordinal between
$\gamma^*(\delta,\beta)$ and $\delta$, then for any $m$ the $m$-walk from $\beta$ to $\alpha$ agrees with the $m$-walk
from $\beta$ to $\delta$ until the penultimate step of the latter walk.  This is a familiar pattern of argument invented by Todor{\v{c}}evi{\'c}.

\begin{corollary}
If $\delta<\beta<\lambda$ and $\cf(\delta)>\aleph_0$, then for any $m<\omega$ and $i<\rho_2^m(\delta,\beta)$ we have
\begin{equation*}
\gamma^*(\delta,\beta)<\alpha<\delta\Longrightarrow \beta^m_i(\alpha,\beta)=\beta^m_i(\delta,\beta).
\end{equation*}
\end{corollary}
\begin{proof}
The proof is a straightforward induction on $i<\rho^m_2(\delta,\beta)-1$ using the definition $\beta^m_i(\alpha,\beta)$.
\end{proof}

The next lemma takes advantage of the fact that for each $\alpha$, $\langle e^m_\alpha:m<\omega\rangle$ forms an increasing family of sets.

\begin{lemma}
If $\alpha<\beta<\lambda$, then there is an $m(\alpha,\beta)<\omega$ such that
\begin{equation*}
\rho^m_2(\alpha,\beta)=\rho^{m(\alpha,\beta)}_2(\alpha,\beta)
\end{equation*}
for all $m\geq m(\alpha,\beta)$. Moreover,
\begin{equation*}
\beta^m_i(\alpha,\beta)=\beta^{m(\alpha,\beta)}_i(\alpha,\beta)
\end{equation*}
for all $i\leq \rho_2^{m(\alpha,\beta)}(\alpha,\beta)$.
\end{lemma}
\begin{proof}
This follows easily by induction once we note that for any ordinals $\alpha<\beta$, the sequence
$\langle \min(e^m_\beta\setminus\alpha):m<\omega\rangle$ is non-increasing and hence eventually constant.
\end{proof}

Armed with the last lemma, we are in a position to define another parameter needed for our construction.

\begin{definition}
Given $\alpha<\beta<\lambda$,   define
 \begin{equation}
 \gamma(\alpha,\beta):= \beta^{m(\alpha,\beta)}_{\rho^{m(\alpha,\beta)}_2(\alpha,\beta)-1}.
 \end{equation}
 \end{definition}

 In plain language, we know that the various $m$-walks from $\beta$ down to $\alpha$ are in agreement for all $m\geq m(\alpha,\beta)$, and so $\gamma(\alpha,\beta)$ is simply the last ordinal visited by all of these walks before they arrive at their destination~$\alpha$.  The last two parameters we need for the proof of Theorem~\ref{filterthm} are given in the following definition; note that we restrict the definition to the case where $\delta\in S$.

\begin{definition}
Suppose $\delta\in S$ and $\delta<\beta<\lambda$. We define

\begin{equation}
\epsilon^*(\delta,\beta):=\text{ least } \epsilon<\sigma\text{ for which }
\gamma^*(\delta,\beta)\leq c^0_\delta(\omega\cdot\epsilon)
 \end{equation}
 and
 \begin{equation}
 m^*(\delta,\beta):=\text{ least $m\geq m(\delta,\beta)$ such that $\cf(\gamma(\delta,\beta))< \mu_m$},
 \end{equation}
 where $c^0_\delta$ and $\langle \mu_m:m<\omega\rangle$ are as in our list of assumptions.
 \end{definition}

The following proposition captures a few facts about the how the various objects we have been considering interact; part~(3) of the proposition has particular importance for us.

 \begin{proposition}
 \label{prop6.4}
 Suppose $\delta\in S$ and $\delta<\beta<\lambda$.
 \begin{enumerate}
 \item $C_\delta[m]\subs e^m_{\gamma(\delta,\beta)}$ for all $m\geq m^*(\delta,\beta)$.
 \sk
 \item Assume $\epsilon^*(\delta,\beta)\leq\epsilon<\sigma$ and $m^*(\delta,\beta)\leq m<\omega$.
 If $\beta^*\in\nacc(C_\delta)\cap I(\epsilon, m)$, then
     \begin{equation*}
     \beta^*\in\nacc(e^m_{\gamma(\delta,\beta)}),
     \end{equation*}
     and
     \begin{equation*}
     \gamma^*(\delta,\beta)\leq\sup(e^m_{\gamma(\delta,\beta)}\cap\beta^*)
     \end{equation*}
 \sk
 \item Furthermore, in  the situation of (2), if $\sup(\beta^*\cap e^m_{\gamma(\delta,\beta)})<\alpha<\beta^*$,
 then
 \begin{equation}
 \beta^m_k(\alpha,\beta)=\beta^m_k(\delta,\beta)\text{ for all }k<\rho^m_2(\delta,\beta),
 \end{equation}
 and
 \begin{equation} \beta^m_{\rho_2^m(\delta,\beta)}(\alpha,\beta)=\beta^*.
 \end{equation}
 In particular, the $m$-walk from $\beta$ to~$\beta^*$ (including the final step!) is an initial segment of the $m$-walk from $\beta$ to~$\alpha$.
 \end{enumerate}
 \end{proposition}
\begin{proof}

For~(1), we know $\delta\in e^m_{\gamma(\delta,\beta)}$ by our definition of $\gamma(\delta,\beta)$ and $m^*(\delta,\beta)$, and a glance back at Lemma~\ref{swallowlemma} gives us what we need.

Condition~(2) gives us some more information --- it claims that some of the structure of $C_\delta$ survives in $e^m_{\gamma(\delta,\beta)}$, and it is at this point that we cash in some of our assumptions for the first time (in particular, it will become clear why we work with the somewhat awkwardly defined $m^*(\delta,\beta)$ instead of $m(\delta,\beta)$).  Given $\beta^*$ in $\nacc(C_\delta)\cap I(\epsilon, m)$, we know that $\beta^*\in e^m_{\gamma(\delta,\beta)}$ by~(1). Condition~(1) of Lemma~\ref{swallowlemma} taken with the definition of $m^*(\delta,\beta)$ tells us
\begin{equation*}
\left|e^m_{\gamma(\delta,\beta)}\right|\leq\cf(\gamma(\delta,\beta))+\mu_m^+\leq\mu_m^+.
\end{equation*}
However, the fact that $\beta^*\in I(\epsilon, m)$ tells us that
\begin{equation*}
\cf(\beta^*)>\mu_m^+,
\end{equation*}
and therefore $\beta^*$ cannot be an accumulation point of $e^m_{\gamma(\delta,\beta)}$.
 This, together with the definition of~$\epsilon^*(\delta,\beta)$, implies
\begin{equation*}
\gamma^*(\beta,\delta)\leq\sup(\beta^*\cap e^m_{\gamma(\delta,\beta)})<\beta^*.
\end{equation*}
In particular, there are many $\alpha$ satisfying the hypothesis of~(3).

Given such an $\alpha$, the fact that $\gamma^*(\delta,\beta)<\alpha$ implies
\begin{equation*}
\beta^m_i(\alpha,\beta)=\beta^m_i(\delta,\beta)
\end{equation*}
for all $i<\rho^m_2(\delta,\beta)$, and therefore
\begin{equation*}
\beta^m_{\rho^m_2(\delta,\beta)-1}(\alpha,\beta) = \beta^m_{\rho^m_2(\delta,\beta)-1}(\delta,\beta) = \gamma(\delta,\beta).
\end{equation*}
Given~(2) and our assumptions about $\alpha$, it follows that
\begin{equation*}
\beta^m_{\rho^m_2(\delta,\beta)}(\alpha,\beta) = \min( e^m_{\gamma(\delta,\beta)}\setminus \alpha)=\beta^*,
\end{equation*}
and so~(3) is established.
\end{proof}

The following corollary follows immediately, and provides a useful upgrade of the preceding proposition.
\begin{corollary}
\label{upgrade}
Suppose $\delta\in S$ and $t$ is a finite subset of $(\delta,\lambda)$, and define
\begin{equation*}
\epsilon^*:=\max\{\epsilon(\delta,\xi):\xi\in t\},
\end{equation*}
and
\begin{equation*}
m^*:=\max\{m^*(\delta,\xi):\xi\in t\}.
\end{equation*}
Given $\beta^*\in \nacc(C_\delta)\cap I(\epsilon, m)$, the ordinal
\begin{equation*}
\gamma^*:=\max\{\sup(e^m_{\gamma(\delta,\xi)}\cap\beta^*):\xi\in t\}<\beta^*
\end{equation*}
has the property that whenever $\gamma^*<\alpha<\beta^*$, we have
\begin{equation}
\beta^m_k(\alpha, \xi) = \beta^m_i(\delta,\xi)\text{ for all $\xi\in t$ and $k<\rho^m_2(\delta,\xi)$,}
\end{equation}
and
\begin{equation}
\beta^m_{\rho^m_2(\delta,\xi)}(\alpha,\beta) = \beta^*\text{ for all $\xi\in t$}.
\end{equation}
In particular, for any $\xi\in t$ the $m$-walk from $\xi$ to $\beta^*$ is an initial segment of the $m$-walk from $\xi$ to $\alpha$.
\end{corollary}
\begin{proof}
This is an immediate consequence of Proposition~\ref{prop6.4}.
\end{proof}

Let us turn now to the proof of Theorem~\ref{filterthm}, and assume $\langle t_\alpha:\alpha<\lambda\rangle$
is a pairwise disjoint family of finite subsets of~$\lambda$.  Define
\begin{equation*}
E:=\{\delta<\lambda: \alpha<\delta\Longrightarrow t_\alpha\subs\delta\}.
\end{equation*}
Since $E$ is closed and unbounded in $\lambda$, we can apply Proposition~\ref{propguess} to conclude that
\begin{equation}
\label{tdefn}
T:=\{\delta\in S: E\cap C_\delta\notin I_\delta\}
\end{equation}
is stationary.

Fix $\delta\in T$, and choose $\beta>\delta$ so that $\delta<\min(t_\beta)$. Just as in Corollary~\ref{upgrade},
we define

\begin{equation*}
\epsilon^*:=\max\{\epsilon(\delta,\xi):\xi\in t_\beta\}
\end{equation*}
and
\begin{equation*}
m^*:=\max\{m^*(\delta,\xi):\xi\in t_\beta\}.
\end{equation*}

The following claim is where the closed unbounded set $E$ starts to become relevant.

\begin{claim}
\label{keyclaim}
Given $\epsilon^*\leq\epsilon<\sigma$,  $m^*\leq m<\omega$, and $\beta^*\in E\cap\nacc(C_\delta)\cap I(\epsilon, m)$, there is an $\alpha$ such that
\begin{equation}
\beta^m_{\rho^m_2(\delta,\xi)}(\zeta,\xi)=\beta^*
\end{equation}
for all $\zeta\in t_\alpha$ and $\xi\in t_\beta$.
\end{claim}
\begin{proof}
Given such $\beta^*$, we work as in Corollary~\ref{upgrade} and define
\begin{equation*}
\gamma^*:=\max\{\sup(e^m_{\gamma(\delta,\xi)}\cap\beta^*):\xi\in t_\beta\}.
\end{equation*}

Since $\beta^*$ is in $E$, we can choose an $\alpha<\lambda$ so that $t_\alpha$ is contained in the interval $(\gamma^*,\beta^*)$.  From part~(3) of Proposition~\ref{prop6.4}, we conclude
\begin{equation*}
\beta^m_{\rho^m_2(\delta,\xi)}(\zeta,\xi)=\beta^*
\end{equation*}
for all $\zeta\in t_\alpha$ and $\xi\in t_\beta$, as required.
\end{proof}

The proof of Theorem~\ref{filterthm} is almost complete.  Given $\delta\in T$, our work provides us
with $\epsilon^*<\sigma$ and $m^*<\omega$ as in Claim~\ref{keyclaim}.  Let us define
\begin{equation*}
B_\delta:= \nacc(C_\delta)\cap \bigcup\{I(\epsilon, m):\epsilon^*\leq \epsilon<\sigma, m^*\leq m<\omega\}.
\end{equation*}
Proposition~\ref{almostallprop} tells us that almost all (in the sense of $I_\delta$) members
 of $C_\delta$ are in $B_\delta$, while Claim~\ref{keyclaim} establishes that for any $\beta^*\in E\cap B_\delta$,
 there are $\alpha<\beta<\lambda$ and $m<\omega$ such that the $m$-walk from $\xi$ to $\zeta$ passes through~$\beta^*$
 for all $\zeta\in t_\alpha$ and $\xi\in t_\beta$.  Thus, if we define
 \begin{equation*}
 A:=\lambda\setminus \bigcup\{E\cap B_\delta:  \delta\in T\},
 \end{equation*}
  the proof Theorem~\ref{filterthm} will be finished if we can show
  \begin{equation*}
  A\in \id_p(\bar{C},\bar{I}).
  \end{equation*}

To establish this, we show
\begin{equation}
\label{eqn12}
A\cap E\cap C_\delta\in I_\delta\text{ for all $\delta\in S$}.
\end{equation}
Clearly we need only worry about those $\delta$ for which $E\cap C_\delta\notin I_\delta$.  Given such a $\delta$, we have
\begin{equation}
\label{eqn13}
A\cap E\cap B_\delta =\emptyset
\end{equation}
by the definition of $A$. If $A\cap E\cap C_\delta$ were $I_\delta$-positive, then this set would have non-empty
intersection with $B_\delta$, as the complement of $B_\delta$ is in $I_\delta$.  This cannot be the case, as we would contradict~(\ref{eqn13}).  Thus, $A\in\id_p(\bar{C},\bar{I})$ and the proof of Theorem~\ref{filterthm} is complete.

\section{Scales and elementary submodels}

One of our assumptions in the last two sections is that we have at our disposal a fixed increasing sequence of regular cardinals $\langle \mu_n:n<\omega\rangle$ that is cofinal in $\mu$.  This parameter was utilized in the definition
of~$m^*$, but otherwise it has been in the background. In the upcoming work, we are going to need to assume that this
sequence of regular cardinals carries {\em a scale}, and since our we are going to need to use quite a bit of scale combinatorics, it seems reasonable to devote some time toward fixing notation and reviewing what we need.

\begin{definition}
Let $\mu$ be a singular cardinal. A {\em scale for
$\mu$} is a pair $(\vec{\mu},\vec{f})$ satisfying
\begin{enumerate}
\item $\vec{\mu}=\langle\mu_i:i<\cf(\mu)\rangle$ is an increasing sequence of regular cardinals
such that $\sup_{i<\cf(\mu)}\mu_i=\mu$ and $\cf(\mu)<\mu_0$.
\item $\vec{f}=\langle f_\alpha:\alpha<\mu^+\rangle$ is a sequence of functions such that
\begin{enumerate}
\item $f_\alpha\in\prod_{i<\cf(\mu)}\mu_i$.
\item If $\gamma<\delta<\mu^+$ then $f_\gamma<^* f_\delta$, where  the notation $f<^* g$  means that $\{i<\cf(\mu): g(i)\leq f(i)\}$ is bounded in $\cf(\mu)$.
\item If $f\in\prod_{i<\cf(\mu)}\mu_i$ then there is an $\alpha<\mu^+$ such that $f<^* f_\alpha$.
\end{enumerate}
\end{enumerate}
Given an increasing sequence $\vec{\mu}$ as above, if there is a $\vec{f}$ such that $(\vec{\mu},\vec{f})$ is a scale,
then say that $\vec{\mu}$ {\em carries a scale}.
\end{definition}

Since we have been working in the situation where the singular cardinal~$\mu$ has countable cofinality, we will deal only with that special case in this section; we refer the reader to~\cite{nsbpr} or~\cite{myhandbook} for more information on these matters.  Thus, from now on we assume

\begin{itemize}
\item $\langle \mu_n:n<\omega\rangle$ is an increasing sequence of regular cardinals such that
\begin{itemize}
\sk
\item $\langle \mu_n:n<\omega\rangle$ is cofinal in $\mu$,
\sk
\item $\sigma<\mu_0$, and
\sk
\item $\langle \mu_n:n<\omega\rangle$ carries a scale $\langle f_\alpha:\alpha<\lambda\rangle$.
\end{itemize}
\end{itemize}
That such sequences exist is a very non-trivial result of Shelah (Chapter~II of~\cite{cardarith}). Interested readers can also consult~\cite{myhandbook} or~\cite{cummings} for expository treatments of this subject.

The function $\Gamma$ defined below is a standard combinatorial tool associated with scales; again, we tailor the definition to the context at hand.

\begin{definition}
For $\alpha<\beta<\lambda$, we define $\Gamma(\alpha,\beta)$ by
\begin{equation}
\label{Gammadef}
\Gamma(\alpha,\beta)=\max\{n<\omega: f_\beta(n)\leq f_\alpha(n)\}.
\end{equation}
For convenience, we say $\Gamma(\alpha,\alpha)=\infty$.
\end{definition}

The function $\Gamma$ has many nice properties ---  it witnesses $\pr_1(\mu^+,\mu^+,\cf(\mu),\cf(\mu))$ (Conclusion~1 on page~67 of~\cite{cardarith}), for example. We shall not make use of this fact directly, but it is certainly in the background throughout our arguments.

The next lemma is a special case of Lemma~7 in~\cite{nsbpr}.  We remind the reader that notation of the form ``$(\exists^*\beta<\lambda)\psi(\beta)$''
  means $\{\beta<\lambda:\psi(\beta)\text{ holds}\}$ is unbounded below~$\lambda$, while
    ``$(\forall^*\beta<\lambda)\psi(\beta)$'' means that $\{\beta<\lambda:\psi(\beta)\text{ fails}\}$ is bounded
     below~$\lambda$.

\begin{lemma}
\label{lemma?}
There is a closed unbounded $C\subs\lambda$ such that the following holds for every $\beta\in C$:
\begin{equation}
\label{eqn15}
(\forall^*n<\omega)(\forall\eta<\mu_n)(\forall\nu<\mu_{n+1})(\exists^*\alpha<\beta)\left[f_\alpha(n)>\eta\wedge f_\alpha(n+1)>\nu\right].
\end{equation}
\end{lemma}

Our conventions regarding elementary submodels are standard. We assume that~$\chi$ is a sufficiently
large regular cardinal and let $\mathfrak{A}$ denote the structure $\langle H(\chi),\in, <_\chi\rangle$
where $H(\chi)$ is the collection of sets hereditarily of cardinality less than $\chi$, and $<_\chi$ is some
well-order of $H(\chi)$.  The use of $<_\chi$ means that our structure $\mathfrak{A}$ has definable Skolem functions,
and we obtain the set of {\em Skolem terms} for $\mathfrak{A}$ by closing the collection of Skolem functions
under composition.

\begin{definition}
Let $B\subs H(\chi)$. Then $\Sk_{\mathfrak{A}}(B)$ denotes the Skolem hull of $B$ in the structure $\mathfrak{A}$.
More precisely,
\begin{equation*}
\Sk_{\mathfrak{A}}(B)=\{t(b_0,\dots,b_n):t\text{ a Skolem term for $\mathfrak{A}$ and }b_0,\dots,b_n\in B\}.
\end{equation*}
\end{definition}

The set $\Sk_{\mathfrak{A}}(B)$ is an elementary substructure of $\mathfrak{A}$, and it is the smallest such structure
containing every element of $B$. The following technical lemma is due originally to Baumgartner~\cite{jb}; it is a
 fact that is quite useful in proving things about the function $\Gamma$ associated with a given scale. Again, we refer
 the reader to~\cite{nsbpr} or~\cite{myhandbook} for a proof.

\begin{lemma}
\label{newcharlem}
Assume that $M\prec\mathfrak{A}$ and let $\sigma\in M$ be a cardinal.  If we define $N=\Sk_{\mathfrak{A}}(M\cup\sigma)$
then for all regular cardinals $\tau\in M$ greater than $\sigma$, we have
\begin{equation*}
\sup(M\cap\tau)=\sup(N\cap\tau).
\end{equation*}
\end{lemma}

As a corollary to the above, we can deduce an important fact about {\em characteristic functions} of models. Once
again, the following definition is but a special case of a more general definition.

\begin{definition}
\label{chardef}
If $M$ is an elementary submodel of $\mathfrak{A}$
such that
\begin{itemize}
\item $|M|<\mu$, and
\sk
\item $\langle \mu_n:n<\omega\rangle\in M$
\end{itemize}
then the {\em characteristic function of $M$} (denoted $\Ch_M$) is the function
with domain~$\omega$ defined by
\begin{equation*}
\Ch_M(n):=
\begin{cases}
\sup(M\cap\mu_n) &\text{if $\sup(M\cap\mu_n)<\mu_n$,}\\
0  &\text{otherwise.}
\end{cases}
\end{equation*}
\end{definition}

In the situation of Definition~\ref{chardef}, it is clear that $\Ch_M$ is an element of the product
 $\prod_{n<\omega}\mu_n$, and furthermore, $\Ch_M(n)=\sup(M\cap\mu_n)$ for all
sufficiently large $n<\omega$.  We can now see that the following corollary follows immediately from Lemma~\ref{newcharlem}.

\begin{corollary}
\label{skolemhulllemma}
Let $M$ be as in Definition~\ref{chardef}.
If $n^*<\omega$ and we define $N$ to be $\Sk_{\mathfrak{A}}(M\cup\mu_{n^*})$,
then
\begin{equation}
\Ch_M\restr [n^*+1,\omega)=\Ch_N\restr [n^*+1,\omega).
\end{equation}
\end{corollary}

We end this section with one more handy bit of terminology due to Shelah~\cite{108}.

\begin{definition}
A $\lambda$-approximating sequence is a continuous $\in$-chain
$\mathfrak{M}=\langle M_i:i<\lambda\rangle$ of elementary submodels of $\mathfrak{A}$ such that
\begin{enumerate}
\item $\lambda\in M_0$,
\item $|M_i|<\lambda$,
\item $\langle M_j:j\leq i\rangle\in M_{i+1}$, and
\item $M_i\cap\lambda$ is a proper initial segment of $\lambda$.
\end{enumerate}
If $x\in H(\chi)$, then we say that $\mathfrak{M}$ is a $\lambda$-approximating sequence over $x$ if
$x\in M_0$.
\end{definition}

Note that if $\mathfrak{M}$ is a $\lambda$-approximating sequence and $\lambda=\mu^+$, then $\mu+1\subs M_0$ because of condition~(4) and the fact that $\mu$ is an element of each $M_i$.

\section{Defining the coloring}

In this section, we will mix ideas from~\cite{535} and~\cite{nsbpr} together with the proof of Theorem~\ref{filterthm}
in order to prove a coloring theorem in {\sf ZFC}.  Our assumptions and notation are as in the previous two sections, and our
goal is to obtain the following theorem.

\begin{theorem}
\label{mainfilterthm}
There is a function $c:[\lambda]^2\rightarrow \lambda$ such that for any pairwise disjoint collection $\langle t_\alpha:\alpha<\lambda\rangle$ of finite subsets of $\lambda$, it is the case that for $\id_p(\bar{C},\bar{I})$-almost
all $\beta^*<\lambda$, we can find $\alpha<\beta$ such that
\begin{equation}
c(\zeta,\xi)=\beta^*\text{ for all }\zeta\in t_\alpha\text{ and }\xi\in t_\beta.
\end{equation}
\end{theorem}

The proof of Theorem~\ref{mainfilterthm} makes use of the function $\Gamma$ defined
 from the scale $\langle f_\alpha:\alpha<\lambda\rangle$. We begin by defining a sequence $\langle c_m:m<\omega\rangle$
 of functions coloring the pairs from~$\lambda$.

\begin{definition}
\label{cmdefn}
Given $\alpha<\beta<\lambda$ and $m<\omega$, we define
\begin{equation*}
c_m(\alpha,\beta)=\beta^m_{k_m(\alpha,\beta)}(\alpha,\beta),
\end{equation*}
where $k_m(\alpha,\beta)$ is the least $k\leq\rho^m_2(\alpha,\beta)$ for which
\begin{equation*}
\Gamma\left(\alpha,\beta^m_k(\alpha,\beta)\right)\neq\Gamma(\alpha,\beta).
\end{equation*}
\end{definition}

In English, the value of of $c_m(\alpha,\beta)$ is the first place on the $m$-walk from $\beta$ down to $\alpha$ where
``$\Gamma$ changes''.  Except for the parameter $m$, this is the same coloring we used in~\cite{nsbpr}.

The following lemma contains the heart of the proof of Theorem~\ref{mainfilterthm}; it shows that the sequence
of colorings $\langle c_m:m<\omega\rangle$ has many of the properties we need.

\begin{lemma}
\label{musst}
If $\langle t_\alpha:\alpha<\lambda\rangle$ is a pairwise disjoint collection of finite subsets of~$\lambda$, then for
$\id_p(\bar{C},\bar{I})$-almost all $\beta^*<\lambda$, there are an $m<\omega$ and $\beta>\beta^*$ such that
\begin{equation}
 \label{eqn8.9}
(\forall^*i<\omega)(\exists^*\alpha<\beta^*)(\forall\zeta\in t_\alpha)(\forall\xi\in t_\beta)\left[\Gamma(\zeta,\xi)=i\wedge c_m(\zeta,\xi)=\beta^*\right].
\end{equation}
\end{lemma}
\begin{proof}
Let $\langle t_\alpha:\alpha<\lambda\rangle$ be given. Clearly we may assume $\alpha\leq\min(t_\alpha)$ for all $\alpha$,
 as we can pass to a subsequence of cardinality $\lambda$ with no loss of generality.  Define $A$ to be the set of $\beta^*<\lambda$ for which
it is {\em impossible} to find an $m<\omega$ and $\beta>\beta^*$ with the required properties. Assume by way of contradiction that
\begin{equation}
\label{contradiction}
A\notin\id_p(\bar{C},\bar{I}).
\end{equation}

 Let $x=\{\bar{C},\bar{e},S,\mu,\lambda, S, \langle \mu_i:i<\omega\rangle, \langle f_\alpha:\alpha<\lambda\rangle,\langle t_\alpha:\alpha<\lambda\rangle\}$ --- all the
 parameters needed to define the sequence $\langle c_m:m<\omega\rangle$ together
  with the sequence of finite sets under consideration --- and let $\mathfrak{M}=\langle M_i:i<\lambda\rangle$
   be a $\lambda$-approximating sequence over $x$.  The set
 \begin{equation*}
 E:=\{\delta<\lambda: M_\delta\cap\lambda=\delta\}
 \end{equation*}
 is closed and unbounded in $\lambda$, so by our assumption~(\ref{contradiction}), there is a $\delta\in S$
 for which
 \begin{equation}
 A\cap E\cap C_\delta\notin I_\delta.
 \end{equation}
 Choose $\beta<\lambda$ so that $\delta<\min(t_\beta)$, and just as in Corollary~\ref{upgrade}, let
 \begin{equation*}
\epsilon^*:=\max\{\epsilon(\delta,\xi):\xi\in t_\beta\}
\end{equation*}
and
\begin{equation*}
m^*:=\max\{m^*(\delta,\xi):\xi\in t_\beta\}.
\end{equation*}
The conclusion of Corollary~\ref{upgrade} tells us that if $\beta^*\in E\cap \nacc(C_\delta)\cap I(\epsilon, m)$ for
some $\epsilon\geq\epsilon^*$ and $m\geq m^*$, then there is an ordinal $\gamma^*<\beta^*$ such that
whenever $t_\alpha$ is contained in the interval $(\gamma^*,\beta^*)$, for every $\zeta\in t_\alpha$ and $\xi\in t_\beta$
we have
 \begin{equation}
 \label{blah}
 \beta^m_{\rho_2^m(\delta,\xi)}(\zeta,\xi)=\beta^*
 \end{equation}
 and
 \begin{equation}
 \label{blah2}
\beta^m_k(\zeta,\xi)=\beta^m_k(\delta,\xi)\text{ for all }k<\rho^m_2(\delta,\xi).
\end{equation}
Proposition~\ref{almostallprop} tells us that $I_\delta$-almost all elements of $C_\delta$ lie in
 \begin{equation*}
B_\delta:= \nacc(C_\delta)\cap \bigcup\{I(\epsilon, m):\epsilon^*\leq\epsilon<\sigma, m^*\leq m<\omega\}
 \end{equation*}
Since $A\cap E\cap C_\delta$ is $I_\delta$-positive,
 it follows that
\begin{equation*}
 A\cap B_\delta\notin I_\delta.
\end{equation*}
 In particular, we can find $\beta^*\in A\cap E\cap C_\delta$ such that $\beta^*\in I(\epsilon, m)$
  for some $\epsilon\geq\epsilon^*$  and $m\geq m^*$, and hence for which is an ordinal $\gamma^*<\beta^*$
with the properties promised in the discussion surrounding~(\ref{blah}) and~(\ref{blah2}).
Our goal is to get a contradiction by establishing that $m$ and $\beta$ witness that this $\beta^*$ is not an element of~$A$.

We observe that since $\beta^*=\sup(M_{\beta^*}\cap\lambda)$, $\beta^*$ must be a member of every closed
unbounded subset of $\lambda$ that is itself an element of $M_{\beta^*}$.  We will shortly make use of this in the
context of Lemma~\ref{lemma?}, but for now we observe the following simple fact:
\begin{equation}
\label{spreadapart}
\alpha<\beta^*\Longrightarrow t_\alpha\subs\beta^*.
\end{equation}

 For each $\alpha<\lambda$, we define a function $f^{\min}_\alpha$ as follows:
 \begin{equation*}
 f^{\min}_\alpha(i)=\min\{f_\zeta(i):\zeta\in t_\alpha\}.
 \end{equation*}
 Since each $t_\alpha$ is finite, it follows that
 \begin{equation*}
 (\forall^*i<\omega)[f^{\min}_\alpha(i)=f_{\min(t_\alpha)}(i)]
 \end{equation*}
 and therefore $\langle f^{\min}_\alpha:\alpha<\lambda\rangle$ is a scale.  Since this new scale is definable from
 parameters available in $M_0$, it itself must be a member of $M_0$ as well.

 As observed earlier, $\beta^*$ is in every closed unbounded subset of $\lambda$ that is also an element of $M_{\beta^*}$.
 Thus, an application of Lemma~\ref{lemma?} inside of $M_{\beta^*}$
 to the scale $\langle f^{\min}_\alpha:\alpha<\lambda\rangle$ tells us~(\ref{eqn15}) holds for~$\beta^*$.
In particular, we can fix $i_0<\omega$ so that whenever $i_0\leq i<\omega$, it is the case that
\begin{equation}
\label{i0}
(\forall \eta<\mu_i)(\forall\nu<\mu_{i+1})(\exists^*\alpha<\beta^*)[f^{\min}_\alpha(i)>\eta\wedge f^{\min}_\alpha(i+1)>\nu].
\end{equation}

The next part of our argument is going to require some Skolem hull arguments. We start by defining
\begin{equation*}
M:=\Sk_{\mathfrak{A}}(x\cup\{\beta^*\})
\end{equation*}
Note that $x\cup\{\beta^*\}\in M_{\beta^*+1}$ (as everything except $\beta^*$ is already in $M_0$) and therefore $M\in M_\delta$
as it  is definable in $M_\delta$ by taking the Skolem hull of $x\cup\{\beta^*\}$ inside the structure $M_{\beta^*+1}$.
Since $M$ is countable and $\mu_0$ is not,
it follows that
\begin{equation*}
\Ch_M(i)=\sup(M\cap\mu_i)<\mu_i\text{ for all }i<\omega,
\end{equation*}
and therefore
\begin{equation}
\label{bounding}
f(i)\leq \Ch_M(i)\text{ for all $f\in M\cap\prod_{n<\omega}\mu_n$ and $i<\omega$.}
\end{equation}

Since $\delta\in E$, it is immediate that $\Ch_M<^* f_\delta$, and the definition of scale tells us that
$\Ch_M<^* f_\gamma$ whenever $\delta\leq\gamma<\lambda$ as well.  Since $t_\beta$ is finite, there is an
$i_1<\omega$ such that
\begin{equation}
\Ch_M\restr [i_1,\omega) < f_{\beta^m_k(\delta,\xi)}\restr [i_1,\omega)\text{ for all $\xi\in t_\beta$ and $k\leq\rho^m_2(\delta,\xi)$.}
\end{equation}
Finally, choose $i_2<\omega$ such that $\cf(\beta^*)<\mu_{i_2}$, and define
\begin{equation*}
i^*=\max\{i_0, i_1, i_2\}.
\end{equation*}
We claim that if $i^*\leq i<\omega$, then
\begin{equation}
(\exists^*\alpha<\beta^*)(\forall \zeta\in t_\alpha)(\forall \xi\in t_\beta)[c_m(\zeta,\xi)=\beta^*\wedge \Gamma(\zeta,\xi)=i].
\end{equation}

\bigskip

Given such an $i$, let $N:=\Sk_{\mathfrak{A}}(M\cup \mu_i)$.
An application of Corollary~\ref{skolemhulllemma} yields
\begin{equation*}
\Ch_N\restr[i+1,\omega)=\Ch_M\restr[i+1, \omega),
\end{equation*}
and since $i_1\geq i$, we conclude from~(\ref{bounding}) that for any $\zeta\in N\cap\lambda$,
\begin{equation}
\label{bounding2}
f_\zeta\restr [i+1,\omega)<f_{\beta^m_k(\delta,\xi)}\restr [i+1, \omega) \text{ for all $\xi\in t_\beta$ and $k<\rho^m_2(\delta,\xi)$.}
\end{equation}

Now define
\begin{equation*}
\eta^*:=\max\{f_{\beta^m_k(\delta,\xi)}(i):\xi\in t_\beta\text{ and }k<\rho^m_2(\delta,\xi)\},
\end{equation*}
and
\begin{equation*}
\nu^*= f_{\beta^*}(i+1).
\end{equation*}
It is clear that $\eta^*<\mu_i$ and $\nu^*<\mu_{i+1}$. Since both these ordinals are, along with $\beta^*$ and the
scale $\langle f^{\min}_\alpha:\alpha<\lambda\rangle$,  elements of $N$, we can apply~(\ref{i0}) inside~$N$
and conclude
\begin{equation}
\label{unbounded}
N\models (\exists^*\alpha<\beta^*)[f^{\min}_\alpha(i)>\eta^*\wedge f^{\min}_\alpha(i+1)>\nu^*].
\end{equation}
Our choice of $i_2$ guarantees $\cf(\beta^*)\subs N$, and therefore $N\cap\beta^*$ is unbounded in $\beta^*$; when
we combine this observation with~(\ref{unbounded}), we conclude
\begin{equation}
\label{realunbounded}
(\exists^*\alpha<\beta^*)[\alpha\in N\wedge f^{\min}_\alpha(i)>\eta^*\wedge f^{\min}_\alpha(i+1)>\nu^*].
\end{equation}

The next proposition will essentially finish our proof.  Recall that ``$\gamma^*$''refers to the ordinal
below $\beta^*$ isolated in the discussion preceding~(\ref{blah}) and~(\ref{blah2}).

\begin{proposition}
\label{prop8.3}
Suppose $\alpha<\beta^*$ satisfies
\begin{itemize}
\item $\alpha\in N$,
\sk
\item $t_\alpha\subs (\gamma^*,\beta^*)$
\sk
\item $f^{\min}_\alpha(i)>\eta^*\wedge f^{\min}_\alpha(i+1)>\nu^*$.
\sk
\end{itemize}
Then for any $\zeta\in t_\alpha$ and $\xi\in t_\beta$, we have
\begin{itemize}
\item $\Gamma(\zeta,\xi)=i$, and
\sk
\item $c_m(\zeta, \xi)=\beta^*$.
\sk
\end{itemize}
\end{proposition}
\begin{proof}
Let $\alpha$ be as hypothesized, and choose $\zeta\in t_\alpha$ and $\xi\in t_\beta$.
We first show
\begin{equation}
\label{gammagoal}
\Gamma(\zeta,\beta^m_k(\zeta,\xi))= i\text{ for all }k<\rho^m_2(\zeta,\xi).
\end{equation}
Given $k<\rho^m_2(\delta,\xi)$, we know that
\begin{equation*}
\beta^m_k(\zeta,\xi)=\beta^m_k(\delta, \xi)
\end{equation*}
by way of~(\ref{blah2}).  Thus,
\begin{equation}
\label{conjunct1}
f_{\beta^m_k(\zeta,\xi)}(i)=f_{\beta^m_k(\delta,\xi)}(i)\leq\eta^*<f^{\min}_\alpha(i)\leq f_\zeta(i).
\end{equation}

Since $t_\alpha\in N$ and $t_\alpha$ is finite, it follows that $\zeta\in N$ as well, and therefore we have
\begin{equation}
\label{conjunct2}
f_\zeta\restr [i+1,\omega)\leq \Ch_M\restr[i+1,\omega)<f_{\beta^m_k(\delta,\xi)}\restr [i+1,\omega)=f_{\beta^m_k(\zeta,\xi)}\restr[ i+1, \omega)
\end{equation}
by way of~(\ref{bounding2}).
The conjunction of~(\ref{conjunct1}) and~(\ref{conjunct2}) tells us $\Gamma(\zeta,\beta^m_k(\zeta,\xi))=i$, and this
establishes~(\ref{gammagoal}).

Since $\beta^m_{\rho_2^m(\zeta,\xi)}=\beta^*$, we are finished if we can establish
\begin{equation}
\Gamma(\zeta,\beta^*)\neq i,
\end{equation}
but this follows easily because we have arranged that
\begin{equation*}
f_\beta^*(i+1)=\nu^*< f^{\min}_\alpha(i+1)\leq f_\zeta(i+1).
\end{equation*}
\end{proof}

We can now finish our proof of~(\ref{eqn8.9}).  Since $\alpha\leq\min(t_\alpha)$, it follows from~(\ref{spreadapart})
that $t_\alpha\subs (\gamma^*,\beta^*)$ whenever $\alpha\in (\gamma^*,\beta^*)$.  Thus, (\ref{realunbounded})
implies that there are unboundedly many $\alpha<\beta^*$ satisfying the assumptions of Proposition~\ref{prop8.3}.
Since $i$ was an arbitrary element of $(i^*,\omega)$, it follows that $m$ and $\beta$ stand witness that
$\beta^*$ is not an element of $A$. This contradiction completes the proof of Lemma~\ref{musst}.
\end{proof}

We turn now to the proof of Theorem~\ref{mainfilterthm}. Our proof requires one more
parameter --- we need to fix a function $p:\omega\rightarrow\omega$ with the property that $p^{-1}(m)$ is infinite
for all $m<\omega$. Given this function, we define our coloring as follows:

\begin{definition}
Given $\alpha<\beta<\lambda$, we define $c:[\lambda]^2\rightarrow\lambda$ by
\begin{equation*}
c(\alpha,\beta) = c_{p(\Gamma(\alpha,\beta))}(\alpha,\beta).
\end{equation*}
\end{definition}

\begin{proof}[Proof of Theorem~\ref{mainfilterthm}]
Suppose $\beta^*$ is as in the conclusion of Lemma~\ref{musst}, and fix $m$ and $\beta$ for which
(\ref{eqn8.9}) is true.  We can find an $i<\omega$ sufficiently large so that
\begin{itemize}
\item $p(i)=m$, and
\sk
\item $(\exists^*\alpha<\beta^*)(\forall\zeta\in t_\alpha)(\forall\xi\in t_\beta)[\Gamma(\zeta,\xi)=i\wedge c_m(\zeta,\xi)=\beta^*]$.
\sk
\end{itemize}
For such $\alpha$, we have
\begin{equation}
c(\zeta,\xi)=\beta^*\text{ for all $\zeta\in t_\alpha$ and $\xi\in t_\beta$},
\end{equation}
just as required by Theorem~\ref{mainfilterthm}.
\end{proof}

We finish this section with a corollary whose proof involves applying a well-known trick
to the coloring from Theorem~\ref{mainfilterthm}.

\begin{corollary}
\label{basiccor}
 Let $\theta\leq\lambda$ be a cardinal. If $\id_p(\bar{C},\bar{I})$ is not weakly $\theta$-saturated, then $\pr_1(\lambda,\lambda,\theta,\aleph_0)$ holds.
 \end{corollary}
 \begin{proof}
Suppose $\pi:\lambda\rightarrow\theta$ is a function such that $\pi^{-1}(\gamma)$ is $\id_p(\bar{C},\bar{I})$-positive for each $\gamma<\theta$.
 Define a coloring $c^*:[\lambda]^2\rightarrow\theta$ by
 \begin{equation*}
 c^*(\alpha,\beta) = \pi(c(\alpha,\beta)).
 \end{equation*}
  Suppose now we are given a pairwise disjoint collection $\langle t_\alpha:\alpha<\lambda\rangle$ of finite subsets of~$\lambda$ and an ordinal $\varsigma<\theta$. Since $\pi^{-1}(\varsigma)$ is $\id_p(\bar{C},\bar{I})$-positive, the conclusion of Theorem~\ref{mainfilterthm} tells us that we can find $\beta^*<\lambda$ and
  $\alpha<\beta<\lambda$ such that
  \begin{itemize}
  \item $\pi(\beta^*)=\varsigma$, and
  \sk
  \item $c\restr t_\alpha\times t_\beta$ is constant with value~$\beta^*$.
  \sk
  \end{itemize}
  Clearly $c^*\restr t_\alpha\times t_\beta$ is constant with value $\varsigma$, and this
  establishes~$\pr_1(\lambda,\lambda,\theta,\aleph_0)$.
 \end{proof}

In light of the above corollary, we see that if $\id_p(\bar{C},\bar{I})$ is not weakly $\mu$-saturated, then $\pr_1(\mu^+,\mu^+,\mu,\aleph_0)$ holds.  A similar situation occurred in~\cite{535} and~\cite{nsbpr}, and in those
two cases we were able to improve things to a coloring with~$\mu^+$ colors and obtain $\pr_1(\mu^+,\mu^+,\mu^+,\cf(\mu))$.  In this paper, however, the use of $m$-walks brings in an extra parameter that seemingly prevents the proofs from the earlier papers from being carried out. We are still able to use a ``stepping-up argument'' to get an upgrade to $\lambda$ colors, but we pay a price in that the resulting coloring is defined on triples instead of pairs.

\begin{theorem}
\label{thm7}
If $\id_p(\bar{C},\bar{I})$ is not weakly $\mu$-saturated, then there exists a coloring $d:[\lambda]^3\rightarrow\lambda$
such that whenever $\langle t_\alpha:\alpha<\lambda\rangle$ is a pairwise disjoint collection of finite subsets of $\lambda$, and $\varsigma<\lambda$, we can find $\alpha<\beta<\gamma$ such that
\begin{equation}
\label{8.20}
(\forall\epsilon\in t_\alpha)(\forall \zeta\in t_\beta)(\forall \xi\in t_\gamma)[ d(\epsilon,\zeta,\xi)=\varsigma].
\end{equation}
\end{theorem}
\begin{proof}
Our main ingredients in the definition of $d$ are the coloring $c$ from Theorem~\ref{mainfilterthm}, and the
 coloring $c^*$ from Corollary~\ref{basiccor}. We also fix functions $g_\beta$ mapping $\beta$ onto $\mu$ whenever $\mu\leq\beta<\lambda$. The function $d:[\lambda]^3\rightarrow\lambda$ is defined by
\begin{equation}
d(\epsilon,\zeta,\xi) =
\begin{cases}
g_{c(\zeta,\xi)}\left(c^*(\epsilon,\xi)\right) &\text{if $\mu\leq c(\zeta,\xi)$, and}\\
0 &\text{otherwise.}
\end{cases}
\end{equation}

Suppose now that $\langle t_\alpha:\alpha<\lambda\rangle$ is a pairwise disjoint family of finite subsets of~$\lambda$. Without loss of generality, we assume $\alpha\leq\min(t_\alpha)$ and $\max(t_\alpha)<\min(t_\beta)$ whenever $\alpha<\beta$; this can easily be arranged by passing to a subfamily of size~$\lambda$, and such a move does not interfere with our conclusion.

Lemma~\ref{musst} and the proof of Theorem~\ref{mainfilterthm} tells us that for $\id_p(\bar{C},\bar{I})$-almost all $\eta<\lambda$, there is value $h(\eta)>\eta$ such that
\begin{equation}
\label{etasatisfies}
(\exists^*\beta<\eta)(\forall\zeta\in t_\beta)(\forall\xi\in t_{h^*(\eta)})\left[ c(\zeta,\xi)=\eta\right].
\end{equation}
Since $\id_p(\bar{C},\bar{I})$ contains all the non-stationary subsets of $\lambda$, Fodor's Lemma implies that there is a single $\iota<\mu$ such that for stationarily many $\eta<\lambda$, we have both~(\ref{etasatisfies}) and that
$g_{\eta}(\iota)$ is defined and equal to $\varsigma$.

Let $T$ denote the set of all $\eta<\lambda$ satisfying the above. By tossing away a non-stationary subset of $T$ if necessary, we can assume that $T$ consists of limit ordinals, and that $t_{h(\eta^*)}\subs \eta$ whenever $\eta^*<\eta$ in $T$. Given this,
if we define
\begin{equation}
s_{\eta}=\{\eta\}\cup t_{h(\eta)},
\end{equation}
then the resulting family $\langle s_\eta:\eta\in T\rangle$ is a pairwise disjoint collection of finite subsets of~$\lambda$.

Since $c^*$ witnesses $\pr_1(\lambda,\lambda,\mu,\aleph_0)$, we can find $\eta^*<\eta$ in $T$ such that $c^*$ is
constant with value $\iota$ when restricted to $s_{\eta^*}\times s_\eta$.
If we define $\alpha = h (\eta^*)$ and $\gamma= h(\eta)$, then we achieve
  \begin{equation}
  \label{8.24}
  (\forall\epsilon\in t_\alpha)(\forall\xi\in t_\gamma)[c^*(\epsilon,\xi)=\iota].
  \end{equation}

Our assumptions on $T$ imply that $\eta = \min(s_\eta)$ and $\max(s_{h(\eta^*)})<\eta$.  Since $\eta$ satisfies (\ref{etasatisfies}), we can choose $\beta<\eta$ such that
\begin{equation}
\max(t_\alpha)<\beta<\eta,
\end{equation}
and
\begin{equation}
\label{8.26}
(\forall\zeta\in t_\beta)(\forall\xi\in t_\gamma)[c(\zeta,\xi)=\eta].
\end{equation}

It is clear that $\alpha<\beta<\gamma$, so we need only verify~(\ref{8.20}).  Suppose then that $\epsilon\in t_\alpha$,
$\zeta\in t_\beta$, and $\xi\in t_\gamma$.  Then $c^*(\epsilon,\xi)=\iota$ and $c(\zeta,\xi)=\eta$ by~(\ref{8.26}) and~\ref{8.24}.  The definition of $T$ implies that $g_\eta(\iota)=\varsigma$, and so
\begin{equation}
d(\epsilon,\zeta,\xi)=g_{c(\zeta,\xi)}(c^*(\epsilon,\xi))=g_\eta(\iota)=\varsigma,
\end{equation}
as required.

\end{proof}

\section{Conclusions}

In this final section, we give a proof of our main theorem. Our goal is to combine Theorem~\ref{thm2} with the coloring theorems from the preceding section (in the case where our singular cardinal has countable cofinality) and from~\cite{nsbpr} (in the case where the cofinality is uncountable). We will dispense with the assumptions that have been in force for the past few sections, in order to state things in full generality. We begin with a short summary of the main results of~\cite{nsbpr}.

\begin{theorem}
\label{oldthm}
Let $\mu$ be a singular cardinal of uncountable cofinality, and let $S$ be a stationary subset of $S^{\mu^+}_{\cf(\mu)}$. There are an $S$-club system $\langle C_\delta:\delta\in S\rangle$ and a coloring $c:[\mu^+]^2\rightarrow\mu^+$ such that, letting $\id_p(\bar{C},\bar{I})$ be defined as in Definition~\ref{idpcidef}, the following hold:
\begin{enumerate}
\sk
\item The ideal $\id_p(\bar{C},\bar{I})$ is $\cf(\mu)$-complete, and $\theta$-indecomposable for every regular cardinal $\theta$ in the interval $(\cf(\mu),\mu)$. (Part~1 of Observation~3.2 on page 139 of~\cite{cardarith})
\sk
\item Whenever $\langle t_\alpha:\alpha<\mu^+\rangle$ is a pairwise disjoint collection of members of $[\mu^+]^{<\cf(\mu)}$, for $\id_p(\bar{C},\bar{I})$-almost all $\beta^*<\mu^+$, there are $\alpha<\beta$ such that $c\restr t_\alpha\times t_\beta$ is constant with value $\beta^*$. (Theorem~2 of~\cite{nsbpr})
\sk
\item If $\id_p(\bar{C},\bar{I})$ is not weakly $\theta$-saturated for $\theta\leq\mu^+$, then $\pr_1(\mu^+,\mu^+,\theta,\cf(\mu))$ holds. (Corollary~19 of~\cite{nsbpr})
\sk
\item If $\pr_1(\mu^+,\mu^+,\mu^+,\cf(\mu))$ fails, then there is an $\id_p(\bar{C},\bar{I})$-positive set $A$ such that the ideal $I$ obtained by restricting $\id_p(\bar{C},\bar{I})$ to $A$ is weakly $\theta$-saturated for some $\theta<\mu$. (Lemma~23 of~\cite{nsbpr})

\sk
\end{enumerate}
\end{theorem}

Notice that our Theorem~\ref{mainfilterthm} gives us the second conclusion of the above theorem in the case where $\cf(\mu)=\aleph_0$. The proof in the case where $\cf(\mu)>\aleph_0$ is much simpler because we can take advantage of the stronger club-guessing theorems known for that case. We will return to the contrast between results from this paper and those from~\cite{nsbpr} after we state and prove our main theorem.

\begin{theorem}[Main Theorem]
Assume $\mu$ is a singular cardinal.
\begin{enumerate}
\sk
\item If $\refl(<\cf(\mu),S^{\mu^+}_{\geq\theta})$ fails for some $\theta<\mu$, then $\pr_1(\mu^+,\mu^+,\theta,\cf(\mu))$ holds.

    \sk

\item If $\refl(<\!\cf(\mu), S^{\mu^+}_{\geq\theta})$ fails for arbitrarily large $\theta<\mu$, then we obtain both $\pr_1(\mu^+,\mu^+,\mu, \cf(\mu))$  and $\mu^+\nrightarrow[\mu^+]^2_{\mu^+}$.

\sk

\item If $\cf(\mu)>\aleph_0$ and $\refl(<\cf(\mu), S^{\mu^+}_{\geq\theta})$ fails for arbitrarily large $\theta<\mu$, then~(2) can be improved to $\pr_1(\mu^+,\mu^+,\mu^+,\cf(\mu))$.
\sk

\item If $\cf(\mu)=\aleph_0$ and $\refl(<\cf(\mu), S^{\mu^+}_{\geq\theta})$ fails for arbitrarily large $\theta<\mu$, then there is a function $d:[\mu^+]^3\rightarrow\mu^+$ such that whenever $\langle t_\alpha:\alpha<\mu^+\rangle$ is a pairwise disjoint family of finite subsets of $\mu^+$ and $\varsigma<\mu^+$, there are $\alpha<\beta<\gamma$ such that
 \begin{equation*}
(\forall\epsilon\in t_\alpha)(\forall \zeta\in t_\beta)(\forall \xi\in t_\gamma)[ d(\epsilon,\zeta,\xi)=\varsigma].
\end{equation*}
\end{enumerate}
\end{theorem}
\begin{proof}

As far as~(1) is concerned, we note that the theorem is only of interest in the case that $\cf(\mu)<\theta$ ---
$\pr_1(\mu^+,\mu^+,\cf(\mu),\cf(\mu))$ holds for any singular cardinal by a result of Shelah (Conclusion~4.1 on page~67 of~\cite{cardarith}). Thus, we assume $\cf(\mu)<\theta$.

Assume $\refl(<\cf(\mu), S^{\mu^+},_{\geq\theta})$ fails for some regular $\theta$ with $\cf(\mu)<\theta<\mu$. If~$\mu$ has uncountable cofinality, then the ideal~$\id_p(\bar{C},\bar{I})$ mentioned in Theorem~\ref{oldthm} cannot be weakly $\theta$-saturated --- if it were, then  part~(1) of Theorem~\ref{oldthm} taken with Theorem~\ref{thm2} would give us $\refl(<\cf(\mu), S^{\mu^+}_{\geq\theta})$. Our result now follows by part~(3) of Theorem~\ref{oldthm}.

What about the case where $\cf(\mu)=\aleph_0$? In this case, our assumptions give us a finite sequence $\langle S_i:i<n\rangle$ of stationary subsets of $S^{\mu^+}_{\geq\theta}$ which fail to reflect simultaneously. Without loss of generality, we can assume that $S_i\subs S^{\mu^+}_{\tau_i}$ for some regular cardinal $\tau_i$. This allows us to choose a regular cardinal~$\sigma$ such that $\max\{\tau_i:i<n\}<\sigma<\mu$.

Let $S= S^{\mu^+}_\sigma$. An application of Theorem~\ref{offcenter} gives us an $S$-club system $\bar{C}$ for which Theorem~\ref{mainfilterthm} holds. The corresponding ideal $\id_p(\bar{C},\bar{I})$ is $\tau_i$-indecomposable for $i<n$ by Proposition~\ref{indecprop}, and so Theorem~\ref{thm2} lets us conclude that $\id_p(\bar{C},\bar{I})$ cannot be weakly $\theta$-saturated --- it it were, then our sets $S_i$ would be stationary subsets of $S^*(\id_p(\bar{C},\bar{I})$ and thus they would reflect simultaneously. We now obtain $\pr_1(\mu^+,\mu^+,\theta,\cf(\mu))$ by way of Corollary~\ref{basiccor}.

As far as part~(2) of our theorem goes, note that Theorem~\ref{wimpythm} gives us the square-brackets part of the conclusion already. We will prove the remainder by contrapositive, focusing only on the case where $\mu$ is of countable cofinality because part~(3) will take care of the other case.

Thus, assume $\cf(\mu)=\aleph_0$ and $\pr_1(\mu^+,\mu^+,\mu,\cf(\mu))$ fails. Apply Theorem~\ref{offcenter} with $S=S^{\mu^+}_{\aleph_1}$ and obtain an $S$-club system $\bar{C}$ for which Theorem~\ref{mainfilterthm} holds. Corollary~\ref{basiccor} implies that $\id_p(\bar{C},\bar{I})$ must be weakly $\mu$-saturated, and Proposition~\ref{indecprop} tells us that the ideal is $\tau$-indecomposable for every regular $\tau$ lying between $\aleph_1$ and $\mu$.

An elementary argument by contradiction establishes the existence of a regular $\theta<\mu$ (without loss of generality greater than $\aleph_1$) and an $\id_p(\bar{C},\bar{I})$-positive set $A$ such that the ideal $I=\{B\subs\mu^+: A\cap B\in\id_p(\bar{C},\bar{I})\}$ is weakly $\theta$-saturated. This ideal $I$ is also $\tau$-indecomposable for any regular $\tau$ with $\aleph_1<\tau<\mu$ --- it inherits this property from $\id_p(\bar{C},\bar{I})$. Thus, $S^*(I)$ is equal to $S^{\mu^+}_{\geq\theta}$ modulo the non-stationary ideal, and $\refl(<\cf(\mu), S^{\mu^+}_{\geq\theta})$ follows by Theorem~\ref{thm2}.

To finish, we note that parts~(3) and~(4) of our main theorem follow by exactly the same argument --- for part~(3) we take advantage of conclusion~(4) of Theorem~\ref{oldthm}, and for part~(4), we use Theorem~\ref{thm7} in the previous section.
\end{proof}

We conclude our paper with a general discussion of some of the issues raised by the research presented here.  First and foremost, it should be clear that we obtain simultaneous reflection almost by accident --- it seems that a much more important phenomenon  is isolated in Theorem~\ref{mainfilterthm} (and Theorem~2 of \cite{nsbpr}), where we get in {\sf ZFC} a coloring theorem intertwined with an easily describable ideal.  In the presence of square-brackets partition relations, the ideal in question must possess some large cardinal type properties and it is not clear at this point if this is even possible:

\begin{question}
Suppose $\mu$ is a singular cardinal. Is it consistent (relative to large cardinals) that ideals of the form $\id_p(\bar{C},\bar{I})$ can be weakly $\theta$-saturated for some $\theta\leq \mu$?
\end{question}

If $2^{\mu}=\mu^+$, then an old theorem of Erd\v{o}s, Hajnal, and Rado~\cite{ehr} tells us that $\mu^+\rightarrow[\mu^+]^2_{\mu^+}$ holds. This motivates the following question, which we phrase in very specific terms:

\begin{question}
Suppose $\mu$ is singular strong limit cardinal of countable cofinality.  Is it consistent (relative to large cardinals) that $2^{\mu}>\mu^+$ and there is a uniform ultrafilter $U$ on $\mu^+$ that is $\theta$-indecomposable for all uncountable regular $\theta<\mu$?
\end{question}

The preceding question may be quite tractable --- for example, if we ignore the cardinal arithmetic aspect of the question, then   Ben-David and Magidor~\cite{indecomposable} show that such a filter can exist on $\aleph_\omega$. We are not sure of the extent to which the need for $2^{\mu}$ to be greater than $\mu^+$ complicates things.

Another somewhat bothersome point is the difference between parts~(3) and~(4) of our main theorem. This discrepancy may be resolvable by an easy argument, but as of yet we do not see how to do it. The general question is as follows:

\begin{question}
Suppose $\mu$ is a singular cardinal for which $\pr_1(\mu^+,\mu^+,\mu,\cf(\mu))$ holds. Does $\pr_1(\mu^+,\mu^+,\mu^+,\cf(\mu))$ hold?
\end{question}

Finally, having established that square-brackets-type relations at successors of singular cardinals necessarily entail simultaneous reflection of stationary sets, we would like to have information about the cofinalities of the ordinals where the reflection takes place:

\begin{question}
Suppose $\mu^+\rightarrow[\mu^+]^2_{\mu^+}$ for $\mu$ singular, and let $\theta<\mu$ be a regular cardinal for which $\refl(<\cf(\mu), S^{\mu^+}_{\geq\theta})$ holds. Can we say anything interesting about the cofinalities of ordinals where the simultaneous reflection takes place?
\end{question}

The answer to the above depends on getting information about the cofinalities of ordinals in the range of $f^*$, where $f^*$ is as in Theorem~\ref{thm2}.

\end{document}